\documentclass[12pt]{amsart}
\usepackage{amsmath,amsfonts,amssymb,amscd,amsthm,amsbsy,epsf,graphicx}

\makeatletter
\@addtoreset{equation}{section}

\makeatother

\newtheorem{theorem}{Theorem}[section]
\newtheorem{lemma}[theorem]{Lemma}

\newtheorem{proposition}[theorem]{Proposition}
\theoremstyle{definition}
\newtheorem*{remark}{Remark}		

\def\myskip{\noalign{\vskip6pt}}

\def\Sp{{S}}
\def\eps{\varepsilon}

\newcommand{\R}{{\mathbb R}}

\newcommand{\lprime}[1]{\,'\!{#1}}
\newcommand{\di}{N}
\newcommand{\twos}{{2^*}}

\newcommand{\ir}{\int_{\R^\di}}
\newcommand{\is}{\int_{S^{\di-1}}}

\newcommand{\jap}[1]{\langle{#1}\rangle}
\newcommand{\rbeta}{r}
\newcommand{\sgn}{\text{sgn}}

\newcommand{\lchoose}[2]{{{#1}\choose{#2}}} 

\begin{document}
\title{On the Boltzmann equation for diffusively excited
granular media}
\author{I. M. Gamba, V. Panferov and C. Villani}
\address{\(^*\)Department of Mathematics,
The University of Texas at Austin,
Austin, TX 78712-1082, USA}
\address{\(^{\dagger}\) UMPA, ENS Lyon, 46 all\'ee d'Italie, 
F-69364 Lyon Cedex 07, France}

\begin{abstract}
We study the Boltzmann equation for a space-homogeneous 
gas of inelastic hard spheres, with a diffusive term 
representing a random background forcing.  
Under the assumption that the initial datum is a 
nonnegative \(L^2(\R^\di)\) function, with bounded mass and 
kinetic energy (second moment), we prove the existence of 
a solution to this
model, which instantaneously becomes smooth and rapidly
decaying. Under a weak additional assumption of bounded third
moment, the solution
is shown to be unique.  We also establish the existence
(but not uniqueness) of a stationary solution.
In addition we show that the high-velocity tails 
of both the stationary and time-dependent particle 
distribution functions are overpopulated with 
respect to the Maxwellian distribution, as conjectured
by previous authors, and we 
prove pointwise lower estimates for the solutions. 
\end{abstract}

\maketitle
\section*{Introduction}
\baselineskip=16pt
In recent years a significant interest has been 
focused on the study of kinetic models for granular 
flows \cite{Ce,Je,Go}. Depending on the external 
conditions (geometry, gravity, interactions with 
surface of a vessel) granular systems may be in a 
variety of regimes, displaying typical features of 
solids, liquids or gases and also producing quite 
surprising effects \cite{UmMeSw}. Finding a systematic 
way to describe such systems under different conditions
is a physical problem of considerable importance. 
At the same time, recent developments in this area 
gave rise to several novel mathematical models with 
interesting properties.  

In the case of rapid, dilute flows, the 
binary collisions between particles may be considered 
the main mechanism of inter-particle interactions 
in the system. In such cases methods of the kinetic 
theory of rarefied gases, based on the Boltzmann-Enskog 
equations have been applied \cite{JeSa,JeRi,GoSh}.

A very important feature of inter-particle 
interactions in granular flows is their
{\it inelastic} character: the total kinetic energy 
is generally not preserved in the collisions. 
Therefore, in order to keep the system out of the 
``freezing'' state, when particles cease to move and the 
system becomes static, a certain driving mechanism, 
supplying the system with energy, is required. 
Physically realistic driven regimes include 
excitation from the moving boundary, through-flow 
of air, fluidized beds, gravity, and other special 
conditions. We accept a simple model for a driving 
mechanism, the so-called thermal bath, in which we 
assume that the particles are subject to uncorrelated 
random accelerations  between the collisions. Such a 
model was studied in \cite{WiMa} in the one-dimensional 
case, and in \cite{VaEr} in general dimension.    

We study the model \cite{VaEr} in the space-homogeneous 
regime, described by the following equation:
\begin{equation}
\label{be:tdep}
\partial_t{f}-\mu\, \Delta_v{f} = Q(f,f),
\quad{v}\in\R^\di,\quad{t>0}.
\end{equation}
Here \(f\) is the one-particle distribution function
(particle density function in the phase space), which is 
a nonnegative function of the microscopic velocity \(v\) 
and the time \(t\); we shall assume $N\geq 2$
(dimension~1 could be treated as well but would require 
a few notational changes).
On the right-hand side of equation \eqref{be:tdep} there is 
the inelastic Boltzmann-Enskog operator for hard spheres 
(the details of which are given below); 
the term \(-\mu\, \Delta_v{f}\), \(\mu=\text{const}\), 
represents the effect of the heat bath. Without loss 
of generality we can set \(\mu=1\) (see Section 
\ref{sim}), which we will from now on assume.
In the sequel, we shall often abbreviate $\Delta_v$ into
just $\Delta$.

One of the interesting features of
the model \eqref{be:tdep} is the fact 
that it possesses nontrivial steady states  described by 
the balance between the collisions and the thermal bath
forcing. Such steady states are given by solutions 
of the equation 
\begin{equation}
\label{be:stat}
\mu\, \Delta_v{f} +  Q(f,f) = 0,\quad v\in\R^\di.
\end{equation}
Solutions of \eqref{be:stat} have been studied in 
\cite{VaEr} by means of formal expansions. The same problem 
was also studied  in \cite{CaCeGa} and in~\cite{BoCe}, 
for a different kind of interactions, namely the 
{\it Maxwell pseudo-particle model}
\cite{BoCaGa,KrBe,KrBe1}, 
by methods of expansions and 
the Fourier transforms, 
respectively. In reference \cite{CeIlSt} the rigorous 
existence of radially symmetric steady solutions for 
the Maxwell model was established.

The aim of this study is to develop a rigorous theory
of for the {\it inelastic hard sphere} model, and to 
investigate the regularity and qualitative properties of the
solutions. 
We prove that equation \eqref{be:tdep} has a unique weak 
solution under basic assumptions that the initial data have 
bounded {\it mass} and {\it kinetic energy}, and satisfy 
some additional conditions (bounded entropy for existence, 
\(L^2(\R^\di)\) for regularity, and bounded third moment 
in \(|v|\) for uniqueness). 
The thermal bath (diffusion) term in \eqref{be:tdep} is 
responsible for the parabolic regularity of solutions: 
the weak solutions become smooth, classical solutions 
after arbitrarily short time.  
We apply generally similar techniques, based on elliptic 
regularity, to treat the steady case. 
Finally, we establish lower bounds, for both steady and  
time-dependent solutions, proving that the distribution
tails are ``overpopulated'' with respect to the Maxwellian,
as was suggested in~\cite{VaEr}. The lower bound for 
steady solutions 
is given by a ``stretched exponential'' 
\(A \exp(-a|v|^{3/2})\), with \(a=a(\alpha,\mu)\). In the
time-dependent case the bound holds with \(A=A(t)\),
where \(A(t)\) is a generally decaying function of time.

We emphasize that the appearance of the ``3/2'' exponent 
is a specific feature of the hard sphere model with 
diffusion, and could 
be predicted by dimensional arguments (cf. \cite{VaEr}).   
On the other hand, the
Maxwell model with diffusion results in a high-velocity 
tail with asymptotic behavior \(C\exp(-c|v|)\), see 
\cite{BoCe}. As a general rule, the exponents in the tails
are expected to depend on the driving and collision mechanisms
\cite{BeCaCaPu, ErBr,ErBr1,BoGaPa}. 
In fact, deviations of the steady states of 
granular systems from Maxwellian equilibria 
(``thickening of tails'') is one of 
the characteristic features of dynamics of granular 
systems, and has been an object of intensive study  
in the recent years \cite{LoCoDeKuGo,KuHe,RoMe,MoShSw}.

We remark that the ``3/2'' bound has rather important
practical implications as well. In particular, it 
indicates that the approximate solutions based on 
the truncated expansion of the deviation from the 
Maxwellian into Sonine polynomials \cite{VaEr,CaCeGa,MoShSw} 
could only be valid for moderate values of \(|v|^2\). 
Any conclusions about the {\it tail behavior} drawn 
from such an expansion should be questioned. 
Indeed, since the deviation function is growing 
rapidly for \(|v|\) large (it is in the weighted 
\(L^1\) space, but not in \(L^2\)\,!),  the Sonine 
polynomial expansion should in general be expected 
to have poor approximation properties in this region. 

The paper is organized as follows. The first section 
contains the preliminaries, where we introduce the 
inelastic collision operator and establish several 
basic identities which are important in the sequel.
In section 2 we establish the bounds for the energy
and entropy of solutions
In Section 3 we study the moments of the distribution 
function by analyzing the {\it moment inequalities}
for equations \eqref{be:tdep} and \eqref{be:stat}. 
The key point in analyzing the moments {is} 
the so-called Povzner inequalities, well-known for 
the classical Boltzmann equation \cite{Po,El,De,We,
Bo,Lu}, which we here extend to the case of inelastic 
interactions and present in a general setting of 
polynomially increasing convex test functions. 
In Section 4 we study the estimates of the inelastic 
collision operator in \(L^p\) spaces with polynomial 
weights, extending the results in \cite{Gu} to the 
inelastic hard sphere case. We continue by establishing 
apriori regularity estimates, based on the interpolation 
of \(L^p\) spaces and the Sobolev-type inequalities. 
In Section 5 we present a rigorous proof of the 
existence and regularity of the time-dependent 
and steady solutions. The arguments presented 
there also justify the formal manipulations 
performed in Sections 2, 3 and 4. In Section 6 we show 
the uniqueness for the time-dependent problem using 
Gronwall's lemma. Finally, in Section 7 we compute 
lower bounds for the stationary and time-dependent 
solutions.

\section{Preliminaries}
\label{sec:prelim}

\subsection{Binary inelastic collisions.}

We study the dynamics of inelastic identical hard balls
with the following law of interactions. 
Let \(v\) and \(v_*\) be the velocities of two particles
before a collision, and denote by \(u=v-v_*\) their 
relative velocity. Let the prime symbol denote the same 
quantities after the collision.  
Then we assume
\begin{equation}
\label{basic-coll}
\begin{split}
&\quad(u'\cdot n )  = - \alpha\, (u\cdot n ),
\\
& u'-(u' \cdot  n) = u -(u\cdot n), 
\end{split}
\end{equation}
where \(n\) is the unit vector in the direction of impact, 
and \(0<\alpha<1\) is a constant called the coefficient 
of normal restitution. Setting \(w=v+v_*\) and using the 
momentum conservation we can express \(v'\) and \(v'_*\) 
as follows:   
\begin{equation}
\label{col:post-1}
v'=\frac{w}{2}+\frac{u'}{2}\,,\qquad
v'_*=\frac{w}{2}-\frac{u'}{2}\,,
\end{equation}
By substituting \eqref{basic-coll} into \eqref{col:post-1}
and equations \eqref{basic-coll}, the post-collisional 
velocities \(v'\) and \(v'_*\) are uniquely determined 
by the pre-collisional ones, \(v\) and \(v_*\),
and the impact parameter \(n\) (cf. \cite{Ce}, \cite{VaEr}). 

The geometry of the inelastic collisions defined by 
relations \eqref{basic-coll}, \eqref{col:post-1} is 
shown in Figure 1. For every \(v\) and \(v_*\) fixed, 
the sets of possible outcomes for post-collisional 
velocities are two (distinct) spheres of diameter 
\(\frac{1+\alpha}{2}|u|\). Thus, it is convenient 
to parametrize the relative velocity after collision 
as follows: 
\begin{equation}
\label{col:post-2}
u'=(1-\beta)\,{u}+\beta\,{|u|}\,\sigma,
\end{equation}
where we denoted \(\beta=\frac{1+\alpha}{2}\). The 
relations \eqref{col:post-1} and \eqref{col:post-2}
define the post-collisional velocities in terms of 
\(v\), \(v_*\) and the angular parameter 
\(\sigma\in S^{\di-1}\).   

\begin{figure}[h] 
\label{figure1}
\begin{center}
\includegraphics[width=12cm]{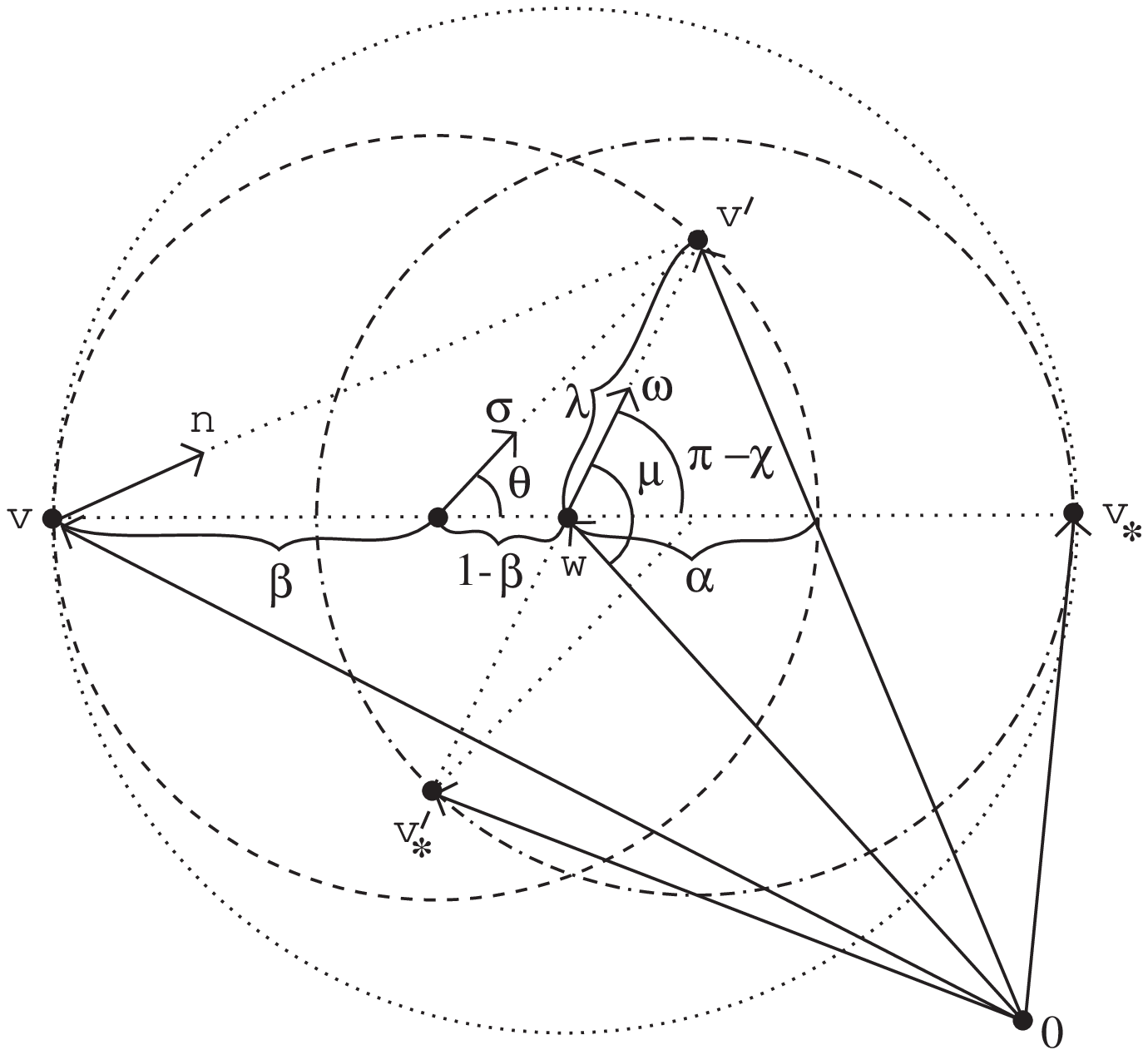}
\end{center}
\caption{
A two-dimensional illustration of the collision mechanism:
\(\scriptstyle{\mathbf -\, -\, -\, -\, -\, -}\,\): possible 
locations of \(v'\); \(\text{}\)
\(- \!\cdot\! - \!\cdot\! -\,\): possible locations of \(v'_*\). 
All lengths shown in assumption \({|u|}/{2} = 1\). 
Unit vectors \(n\), \(\sigma\), \(\omega\) not to scale. } 
\end{figure} 


\subsection{Weak form of the collision operator.}
We define the collision operator by its action on 
test functions, or {\it observables}.
Taking 
\(\psi=\psi(v,t)\) to be a suitably regular test function, 
we introduce the following weak bilinear form of the collision term:
\begin{equation}
\label{co:bil-weak-sym}
\ir Q(g,f)\, \psi\, dv = 
\ir\ir\is f g_* \,(\psi'-\psi)\,|u|\,
b(u,\sigma)\,d\sigma\,dv\,dv_*.
\end{equation}
Here and below we use the shorthand notations \(f=f(v,t)\), 
\(g_*=g(v_*,t)\), \(\psi'=\psi(v',t)\), etc. 
The function \(b(u,\sigma)\) in \eqref{co:weak-sym} 
is the product of the Enskog correlation factor 
\(k(\rho,d)\) (which is a constant in the 
space-homogeneous case) by the 
differential collision cross-section, expressed 
in the variables \(u\), \(\sigma\). In the case of 
hard-sphere interactions,   
\[
b(u,\sigma) = k(\rho,d)\left(\frac{d}{2}\right)^{N-1}
\left(\frac{1-(\nu\cdot\sigma)}{2}\right)^{-\frac{N-3}{2}},
\]
where \(\nu=u/|u|\), and \(d\) is the diameter of the 
particles. Notice that the hard sphere cross-section 
depends only on the angle between \(u\) and \(\sigma\), 
and is generally anisotropic, unless \(N=3\). Without 
restricting generality, by choosing 
the value of \(d\) accordingly, we can always assume 
that
\begin{equation}
\label{norm:xsect}
\is b(u,\sigma)\,d\sigma = 1. 
\end{equation}

{ Of course,
to write down the Boltzmann operator we only need $Q(f,f)$, but
later on it will be sometimes convenient to work with the 
bilinear form $Q(g,f)$. 
An explicit form of $Q$ will be given later on;
however for many purposes it will be easier to work with the
weak formulation which is also quite natural from the
physical point of view (it is analogous to the well-known Maxwell 
form of the Boltzmann collision operator 
\cite[Chapter~1, Section~2.3]{V}).}

In the case when $f=g$ in~\eqref{co:bil-weak-sym}, we can
further symmetrize and write
\begin{equation}
\label{co:weak-sym}
\ir Q(f,f)\, \psi\, dv 
= \frac{1}{2}
\ir\ir\is f f_* \,(\psi'+ \psi'_*-\psi-\psi_*)\,|u|\,
b(u,\sigma)\,d\sigma\,dv\,dv_*.
\end{equation}
Notice that the particular form of the inelastic 
collision laws enters \eqref{co:weak-sym} only 
through the test function \(\psi'\). 

\subsection{Equations for observables and conservation relations.}

Using the weak form \eqref{co:weak-sym} allows us to study 
equations for average values of observables given by 
the functionals of the form \(\ir f \psi\,dv\). Namely, 
multiplying equation \eqref{be:tdep} by a test function 
\(\psi(v,t)\) and integrating by parts we obtain
\begin{equation}
\label{be:weak}
\begin{split}
\Big[ \ir f\,\psi\,dv \Big]_{t=0}^{t=T}
- \int_0^T \!\ir f\,(\partial_t\,\psi 
+ \Delta_v \,\psi) \,dv\,dt
= \int_0^T \!\ir Q(f,f)\, \psi\,dv\,dt. 
\end{split}
\end{equation}
With the weak form  \eqref{co:weak-sym} of the collision 
operator, it is easy to verify formally the basic conservation 
relations that follow from \eqref{be:tdep}. Namely, setting 
\(\psi=1\) and \(\psi=v_i\) in \eqref{be:weak} and assuming 
that \(\ir f\,\psi\,dv\) is differentiable in \(t\), we 
obtain the {\it conservation of mass and momentum}:
\begin{equation}
\label{cons:mm}
\frac{d}{dt}\ir f \{1,v_1,\dots,v_\di\}\,dv = 0 .
\end{equation}
Further, taking \(\psi=|v|^2\) and computing
\begin{equation}
\label{ene-iden}
|v'|^2+|v'_*|^2-|v|^2-|v_*|^2
=-\frac{1-\alpha^2}{2}\,
\frac{1-(\nu\cdot\sigma)}{2}|u|^2,
\end{equation}
we obtain the following relation for the {\it dissipation of 
kinetic energy}:
\begin{equation}
\label{diss:ke}
\frac{d}{dt}\ir f\,|v|^2\,dv = 2\di 
- \epsilon_{\di}\frac{1-\alpha^2}{4}{\ir}{\ir} 
ff_*|u|^3\,dv_*\,dv, 
\end{equation}
where 
\[
\epsilon_\di = \is \frac{1-(\nu\cdot\sigma)}{2}\,
b(u,\sigma)\, d\sigma = \text{const}. 
\]
Notice that, unlike the no-diffusion case, the kinetic 
energy is not necessarily a monotone function of time. However, 
it is not difficult to show using \eqref{diss:ke} 
(see Section \ref{apr:basic}) that the kinetic energy 
remains bounded for all times, provided the initial distribution 
function has finite energy. 

Finally, equation \eqref{be:weak} allows us to define the concept 
of solutions of \eqref{be:tdep} which we use throughout the paper. 
Namely, we say that a function \(f\) 
is a {\it weak solution} of \eqref{be:tdep} if  for every 
\(T>0\), \(f\in L^1([0,T]\times\R^\di)\), 
\(Q(f,f)\in L^1([0,T]\times\R^\di)\) and \eqref{be:weak} holds 
for every \(\psi\in C^1([0,\infty),C^2(\R^\di))\) vanishing 
for \(t>T\). It can be 
shown in the usual way that if a weak solution is 
sufficiently smooth { (say, continuously differentiable with
respect to time and twice continuously differentiable with respect to
velocity)} and satisfies { suitable} decay conditions 
for large \(|v|\), then it also is a classical solution.    

\subsection{Entropy identity}

Taking in the weak form \eqref{co:weak-sym} 
\(\psi=\log f\) we 
obtain an interesting identity for the 
entropy \(\ir f\log f\,dv\). First, we 
compute
\begin{equation}
\label{co:weak-logf}
\begin{split}
{\ir}{Q(f,f)}\log{f}\,dv
=\frac{1}{2}{\ir}{\ir}
{\is}ff_*\log\frac{f'f'_*}{ff_*}\,
|u|\,b(u,\sigma)\,d{\sigma}\,dv\,dv_*
\\
=\frac{1}{2}{\ir}{\ir}
{\is}ff_*\left(\log\frac{f'f'_*}{ff_*}
-\frac{f'f'_*}{ff_*}+1\right)
|u|\,b(u,\sigma)\,d{\sigma}\,dv\,dv_*
\\
+\frac{1}{2}{\ir}{\ir}\is
(f'f_*'-ff_*)\,|u|\,b(u,\sigma)\,d\sigma\,dv\,dv_*.
\end{split}
\end{equation}
The last term vanishes in the elastic case \(\alpha=1\); 
however, as we see below, it is generally different from 
zero if \(\alpha<1\). 
To compare the integral of \(f'f'_*\) to that of \(ff_*\)
we perform the transformation corresponding to 
the {\it inverse collision}, passing from the velocities 
\(v'\), \(v'_*\) to their predecessors \(v\) and \(v_*\). 
Such a transformation is more easily expressed in the 
variables \(u\) and \(n\).  Passing to these variables, 
we can write the integral of 
\(f'f'_*\) as follows:
\begin{equation}
\label{int:prime}
d^{\di-1}{\ir}{\ir}
\int_{S^{\di-1}_+} 
f'f'_*\, | u\cdot n |\,dn\,dv\,dv_*,
\end{equation}
where \(S^{\di-1}_+=\{n\in S^{\di-1}\:|\:u\cdot n>0\}\). 
The ``inverse collision'' transformation 
\((v,v_*,n)\mapsto(v',v'_*,-n)\)
has the Jacobian determinant equal to \(\alpha\) \cite{Ce}. 
Therefore, using the first of the equations \eqref{basic-coll}, 
the integral \eqref{int:prime} is computed as 
\begin{equation}
\label{qplus:alpha}
d^{\di-1}\,\frac{1}{\alpha^2}\ir\ir\int_{S^{\di-1}_+}\,\! 
f f_*  \,
| u\cdot n |\,dn\,dv\,dv_*,
\end{equation}
Changing variables in the angular integral from \(n\) 
to \(\sigma\), we rewrite \eqref{int:prime} as 
\begin{equation}
\label{int:prime-t}
\frac{1}{\alpha^2} \ir\ir\is  f f_*\,
|u|\,b(u,\sigma)\,d\sigma\,dv\,dv_* 
= \frac{1}{\alpha^2} \ir\ir  f f_* |u|\,dv\,dv_* \,.
\end{equation}
In view of \eqref{co:weak-logf} and \eqref{int:prime-t}
the entropy equation becomes
\begin{equation}
\label{diss:entropy}
\begin{split}
\frac{d}{dt}{\ir}f\log{f}\,dv
+  4{\ir} & \left|\nabla\sqrt{f}\right|^2\,dv
\\
= \frac{1}{2}{\ir}{\ir}{\is} & 
ff_*\left(\log\frac{f'f'_*}{ff_*}
-\frac{f'f'_*}{ff_*}+1\right)
|u|\,b(u,\sigma)\,d\sigma\,dv\,dv_*
\\
&\qquad\quad + \,\frac{1}{2}\,\Big(\frac{1}{\alpha^2}-1\Big)\,
{\ir}{\ir}ff_*\,|u|\,dv\,dv_*.
\end{split}
\end{equation}
{In these equations as in all the sequel, the symbol $\nabla$
will stand for the gradient operator with respect to velocity variables.}
Here the first term on the right-hand side is 
nonpositive (notice the inequality \(\log{x}-x+1 \le 0\))
and similar to the entropy dissipation in the elastic case. 
The last term in \eqref{diss:entropy} is a nonnegative 
correction term that vanishes in the elastic limit 
\(\alpha\to1\). 

\subsection{Similarity in the equations 
and normalization of solutions.}
\label{sim}

As a consequence of \eqref{cons:mm},
the total density (mass) and momentum ({mean} value)
of the distribution function are equal to those of the 
initial distribution. We can write this as follows:
\[
\ir f\, dv = \rho_0 = \text{const}, \quad \text{and} 
\quad \ir f\, v_i \, dv = \rho_0 {v_0}_i = \text{const}_i, 
\quad i = 1,\dots,\di . 
\]
In fact, we can always assume that \(\rho_0=1\), \(v_0=0\)
and \(\mu=1\) in \eqref{be:tdep}. Indeed, if \(f(v,t)\) is 
such a solution to \eqref{be:tdep}, then, for every \(\rho_0\), 
\(v_0\) and \(\mu\), the function
\[
f_{\{\rho_0,v_0,\mu\}} (v,t)
= \rho_0 \eta^{-N} f\big(t/\tau,(v-v_0)/\eta\big),
\] 
where 
\[
\tau = \rho_0^{-2/3} \mu^{-1/3}, \quad \text{and} \quad 
\eta = \rho_0^{-1/3} \mu^{1/3},
\]
is a solution corresponding to the given values of \(\rho_0\), 
\(v_0\) and \(\mu\). 

\subsection{Strong form of the collision operator.}

Using the weak form \eqref{co:weak-sym} we can 
derive the usual strong form of the collision operator. 
We notice the obvious splitting into the ``gain'' and 
the ``loss'' terms,
\[
Q(g,f) = Q^+(g,f) - Q^-(g,f). 
\]
Assuming that \(f\) is regular enough, setting 
\(\psi(v)=\delta(v-v_0)\) in the part of 
\eqref{co:weak-sym} corresponding to \(Q^-(g,f)\), 
and using \eqref{norm:xsect} we find
\[
Q^-(g,f) = \ir\is fg_* \,|u|\,b(u,\sigma)\,d\sigma\,dv_* = f\,(g*|v|).
\]

To find the explicit form of \(Q^+(g,f)\) we invoke
the inverse collision transformation, tracing the collision 
history back from the pair \(v\), \(v_*\) to 
their predecessors, which we denote 
by \(\lprime{v}\) and \(\lprime{v_*}\). 
Setting \(\psi(v) = \delta(v-v_0)\) and arguing similarly 
to the derivation of the entropy identity we obtain
\[
Q^+(g,f) = \ir\is \lprime{f}\lprime{g_*}\,\frac{1}{\alpha^2}\,
|u|\,b(u,\sigma)\,d\sigma\,dv_*\,,
\]
where \(\lprime{f}=f(\lprime{v},t)\), \(\lprime{g_*} 
= g(\lprime{v_*},t)\), and the pre-collisional velocities 
are defined as 
\begin{equation}
\label{col:pre}
\lprime{v}=\frac{w}{2}+\frac{\lprime{u}}{2}\,,
\quad
\lprime{v_*}=\frac{w}{2}-\frac{\lprime{u}}{2}\,,
\quad\text{where}\quad
\lprime{u}
=(1-\gamma){u} + \gamma\,{|u|}\sigma,\qquad
\end{equation}
and \(\gamma  = \frac{\alpha+1}{2\alpha}\).

\section{Basic Apriori Estimates: Energy and Entropy}
\label{apr:basic}

{In the classical theory of the elastic Boltzmann equation,
the {\em energy conservation} and the {\em entropy decay} are 
the most fundamental facts which provide the base for every 
analysis.  In the present setting naturally we do not have 
energy conservation, and the energy inequality (expressing 
that collisions do not increase the energy) would by no means 
be sufficient to compensate for that. So the
key ingredient will be to replace it by the more precise
{\em energy dissipation} estimate, as follows.}

To study  solutions of \eqref{be:tdep} and \eqref{be:stat}
we assume for simplicity that they satisfy the normalization 
conditions of unit mass and zero average; however the estimates 
we derive below will be by no means restricted to such solutions. 
We use the  energy equation \eqref{diss:ke} and apply
Jensen's inequality for the last term  to get
\[
{\ir}f_*|u|^3\,dv_* \ge \Big|v-\ir f(t,v)\, v \, dv\Big|^3=|v|^3,
\]
and therefore,
\begin{equation*}
{\ir}{\ir} 
ff_*|u|^3\,dv_*\,dv \ge \ir f |v|^3 \, dv.  
\end{equation*}
We then get (in the time-dependent case) the differential inequality
\begin{equation}
\label{mom:2}
\frac{d}{dt}\ir f\,|v|^2\,dv + k_1 
\ir f\,|v|^3\,dv \le K_1,
\end{equation}
where \(K_1=2\di\) and 
\(k_1 = {\epsilon_\di} \frac{1-\alpha^2}{4}\). Further,
by Jensen's inequality,
\[
\ir f \,|v|^3 \, dv \ge \Big( \ir f \,|v|^2 \, dv \Big)^{3/2} ,
\]
and we obtain 
\[
\frac{d}{dt} \ir f \,|v|^2 \, dv  \le
K_1 - k_1 \Big( \ir f \,|v|^2 \, dv \Big)^{3/2}.
\]
Thus, if \(\ir f |v|^2\,dv > (K_1/k_1)^{2/3}\), then 
\(\frac{d}{dt}\ir f |v|^2\,dv<0\) and so,
\begin{equation*}
\sup\limits_{t\ge0}\ir f\,|v|^2\,dv \le \max \Big\{\ir f_0\,|v|^2\,dv,
\big({K_1}/{k_1}\big)^{{2}/{3}}\Big\}.
\end{equation*}
In the steady case the derivative term drops in \eqref{mom:2},
and we obtain 
\[
\ir f\,|v|^3 dv \le \frac{K_1}{k_1}.  
\] 

Let us introduce
the following weighted \(L^1\) spaces: 
\begin{equation}
\label{L1k}
L^1_k(\R^\di)=\{f\:|\:f\jap{v}^{k}\in{L^1}(\R^\di)\},
\end{equation}
where \(k\ge0\) and \(\jap{v}=(1+|v|^2)^{1/2}\). We then 
define 
the norms in \(L^1_k\) as \(\ir |f| \jap{v}^k\,dv\), which 
for \(f\) nonnegative coincide with the {\it moments} 
\(\ir f \jap{v}^k\,dv\). The above argument implies 
apriori estimates for the steady solutions 
in \(L^1_3(\R^\di)\), and for the time-dependent ones in 
\(L^\infty([0,\infty),L^1_2(\R^\di))\) and 
\(L^1_{\rm loc}([0,\infty),L^1_3(\R^\di))\). 
We emphasize that the bounds depend on \(\alpha\)
and deteriorate in the elastic limit \(\alpha\to1\). 
In fact, these bounds for \(\alpha<1\) make a most 
striking contrast 
with the classical Boltzmann equation for elastic 
particles. 

Next, using the entropy equation 
\eqref{diss:entropy} we show that the entropy is 
bounded uniformly in time, for initial data with
finite mass, kinetic energy and entropy. 
To obtain this, we first estimate the second term in 
\eqref{diss:entropy} using the Sobolev embedding 
inequality: { assuming for simplicity here that
$N\geq 3$, we have}
\[
\ir |\nabla\sqrt{f}|^2\,dv \ge c \|f\|_{L^{p^*}} ,
\]
where \(p^*=\di/(\di-2)\). Further, we have the inequality
\begin{equation}
\label{ineq:logeps}
\ir f\log f \, dv \le C_\eps\|f\|_{L^{p^*}}^{\eps} \, ,
\end{equation}
for all \(\eps>0\). Indeed, obviously,
for every \(\delta>0\),
\begin{equation}
\label{ineq:logdelta}
\ir f\log f \, dv \le C_\delta \ir f^{1+\delta}\,dv.
\end{equation}
Further, by H\"older's inequality, for \(\delta<p^*\),
\[
\|f\|_{L^{1+\delta}} \le \|f\|^{1-\nu}_{L^1} \|f\|^{\nu}_{L^{p^*}},
\]
where \(\nu=\frac{p^*\delta}{(p^*-1)(1+\delta)}\). Therefore, 
\[
\ir f^{1+\delta}\,dv \le \|f\|_{L^{p^*}}^{\frac{p^*}{p^*-1}\delta},
\]
which together with \eqref{ineq:logdelta} implies 
\eqref{ineq:logeps}. Now, coming back to estimating the 
terms in the entropy equation \eqref{diss:entropy}, we get
\begin{equation}
\label{entr:Gron}
\frac{d}{dt}\ir f\log{f}\, dv + c_\eps\Big(\ir f\log f \, dv\Big)^{1/\eps}
\le C \|f\|^2_{L^1_1}\,.
\end{equation}
The established bound in 
\(L^\infty([0,\infty),L^1_2(\R^\di))\) implies that
for initial data with finite mass and energy, the right-hand 
side of \eqref{entr:Gron} is bounded by a constant, and we obtain 
by Gronwall's lemma, 
\[
\sup\limits_{t\ge0}\ir f\log f\, dv 
\le C({\textstyle\ir f_0\log f_0\, dv}, \|f_0\|_{L^1_2})\,. 
\]
Integrating \eqref{diss:entropy} in time, we also get
\(
\sqrt{f}\in L^2([0,T],H^1(\R^\di)),
\) 
for every \(T>0\), which implies in particular, 
\(
f \in L^{p^*}([0,T]\times\R^\di), 
\)
where the constants in the estimates depend on the initial mass, 
energy and entropy of the solutions. For the steady solutions 
we obtain a particularly simple estimate 
\(
\|\nabla \sqrt{f}\|_{L^2} \le C \|f\|_{L^1_1}. 
\)
{As the reader will easily check, our assumption that $N\geq 3$ is
just for convenience, and can easily be circumvented in dimension~2
by the Moser-Trudinger inequality, or just the local control of
all $L^p$ norms of $f$ by $\|\nabla\sqrt{f}\|_{L^2}$, together with
a moment-based localization argument.}

\section{Moment inequalities}

We further look for apriori estimates of the solutions
in the spaces \(L^1_k\) \eqref{L1k} with \(k>2\). 
Such estimates will play a very important role in 
our study of regularity, which we perform in Section 
\ref{sec:reg}. The key technique for 
obtaining the necessary estimates is the  
so-called Povzner inequalities \cite{Po,El,De,MiWe,Bo,Lu} 
which we here extend to the inelastic case. 


\subsection{The Povzner-type inequalities}

We take \(\psi(x)\), \(x>0\) to be a convex 
nondecreasing function and look for estimates 
of the expressions 
\begin{equation}
\label{qpsi}
q\,[\psi](v,v_*,\sigma) = \psi(|v'|^2) + \psi(|v'_*|^2) 
- \psi(|v|^2) - \psi(|v_*|^2)
\end{equation}
and 
\begin{equation}
\label{qpsi_bar}
\bar{q}\,[\psi](v,v_*) = \is (\psi(|v'|^2) + \psi(|v'_*|^2) 
- \psi(|v|^2) - \psi(|v_*|^2))\,b(u,\sigma)\,d\sigma, 
\end{equation}
which appear in the weak form of the 
collision operator \eqref{co:weak-sym}. 

Our aim is to treat the cases of
\begin{equation}\label{pow}
\psi(x)=x^p, \quad \text{and} \quad \psi(x)=(1+x)^p-1, \quad p > 1, 
\end{equation}
and also truncated versions of such functions which will be 
required in the rigorous analysis of moments in 
Section \ref{sec:exist}.  Thus, we will 
require functions \(\psi\) to satisfy the following list 
of conditions:  
\begin{eqnarray}
\label{cond:psi-1} & \psi(x)\ge0, \quad x>0; \quad \psi(0) = 0;
\\ \label{cond:psi-2} & \psi(x) \;\text{is convex,}\; 
C^1([0,\infty)),\; \psi''(x)\; \text{is locally bounded};
\\ \label{cond:psi-3} & \psi'(ax) \le \eta_1(a)\, 
\psi'({x}),\quad x>0,\quad a>1;
\\ & \label{cond:psi-4} \psi''(ax) \le \eta_2(a)\, \psi''(x), 
\quad x>0\quad a>1,
\end{eqnarray}
where \(\eta_1(a)\) and \(\eta_2(a)\) are functions of \(a\) only, 
bounded on every finite interval of  \(a>0\). The above 
conditions are easily verified for the functions 
\eqref{pow}. 

We will further { establish} the following elementary lemma. 

\begin{lemma}
\label{psi:ineq}
Assume that \(\psi(x)\) satisfies 
\eqref{cond:psi-1}--\eqref{cond:psi-4}. 
Then
\begin{equation}
\label{by_above}
\psi(x+y) - \psi(x) - \psi(y) 
\le A\, (\,x\,\psi'(y) + y\,\psi'(x)\,) 
\end{equation}
and
\begin{equation}
\label{by_below}
\psi(x+y) - \psi(x) - \psi(y) \ge b\, xy \,\psi''(x+y). 
\end{equation}
where \(A=\eta_1(2)\) and \(b=(2\eta_2(2))^{-1}\).  
\end{lemma}
\begin{proof} To establish the first of the bounds assume that 
\(x \ge y \). Then, since \(\psi(y) \ge 0\), 
\[
\begin{split}
\psi(x+y) - \psi(x) - \psi(y) 
& \le \psi(x+y) - \psi(x) 
\\
& = \int_0^y \psi'(x+t)\, dt \le \int_0^y \eta_1(2)\,\psi'(x)\, dt 
= A\, y\,\psi'(x), 
\end{split}
\]
By symmetry we have 
\[
\psi(x+y) - \psi(x) - \psi(y)  \le A\, x\,\psi'(y), 
\]
when \(x \le y\).  This proves the required inequality for 
all \(x\) and \(y\).  
To prove the second of the bounds in the lemma, 
we can write, using \eqref{cond:psi-4} and the 
normalization \(\psi(0)=0\), 
\[
\begin{split}
\psi(x+y) - \psi(x) - \psi(y) 
= \int_0^y ( \psi'(x+t) - \psi'(t) ) \, dt
= \int_0^y \int_0^x \psi''(t+\tau) \,d\tau \, dt 
\\
\ge (\eta_2(2))^{-1}\,\psi''(x+y)\int_{0}^y 
\int_{0}^x \chi_{\{t+\tau>(x+y)/2\}} \,d\tau \,dt
= (2\eta_2(2))^{-1}\,xy\, \psi''(x+y). 
\end{split}
\] 
This completes the proof. 
\end{proof}

{In the sequel, we shall use} some relations 
involving post-collisional velocities \(v'\) and \(v'_*\). 
It becomes more convenient to parametrize them in the 
{\it center of mass--relative velocity} variables. We 
therefore set
\begin{equation}
\begin{split}
\label{cm:lambda:1}
v' = \frac{w + \lambda|u|\omega}{2},\quad\text{and}\quad
v'_* = \frac{w - \lambda|u|\omega}{2},
\end{split}
\end{equation}
where \(w=v+v_*\), \(u=v-v_*\), and \(\omega\) 
is a parameter vector
on the sphere \(S^{\di-1}\) (see Figure 1). We have 
\[
\lambda\omega=\beta\sigma+(1-\beta)\nu,
\]
where \(\beta=\frac{1+\alpha}{2}\) and \(\nu=u/|u|\), 
and therefore, 
\begin{equation}
\label{lambda:def}
\lambda=\lambda(\cos\chi)=
(1-\beta)\cos\chi+\sqrt{(1-\beta)^2(\cos^2\chi-1)+\beta^2},
\end{equation}
where \(\chi\) is the angle between \(u\) and \(\omega\). 
Notice that 
\begin{equation*}
0<\alpha \le \lambda(\cos\chi) \le 1,
\end{equation*}
for all \(\chi\). With this parametrization we have 
\begin{equation}
\label{ine-mu}
\begin{split}
|v'|^2&=\frac{|w|^2+\lambda^2|u|^2
+2\lambda|u||w|\,\cos\mu}{4}, \\
|v'_*|^2&=\frac{|w|^2+\lambda^2|u|^2
-2\lambda|u||w|\,\cos\mu}{4},
\end{split}
\end{equation}
where \(\mu\) is the angle between the vectors \(w=v+v_*\) 
and \(\omega\). 

\begin{lemma}
\label{povz:conv}
Assume that the function \(\psi\) 
satisfies \eqref{cond:psi-1}--\eqref{cond:psi-4}. 
Then we have 
\[
q\,[\psi] = - n\,[\psi] + p\,[\psi],
\]
where
\[
p\,[\psi] 
\le A\, ( |v|^2\,\psi'(|v_*|^2) + |v_*|^2\,\psi'(|v|^2)\,) 
\]
and
\[
n\,[\psi] \ge \kappa(\lambda,\mu)\, ( |v|^2 + |v_*|^2 )^2\, 
\psi''(\,|v|^2 + |v_*|^2).
\]
Here \(A\) is the constant in estimate \eqref{by_above}, 
\[
\kappa(\lambda,\mu)
=\frac{b}{4}\,\lambda^2\,(\eta_2(\lambda^{-2}))^{-1}\,
\sin^2\!\mu,
\] 
and \(b\) is the constant in estimate \eqref{by_below}.
\end{lemma}
\begin{proof}
We start by setting 
\[
p\,[\psi] = \psi(|v|^2+|v_*|^2) - \psi(|v|^2) - \psi(|v_*|^2)   
\]
and
\[
n\,[\psi] = \psi(|v|^2+|v_*|^2) - \psi(|v'|^2) - \psi(|v'_*|^2).  
\]
The estimate for \(p\,[\psi]\) follows easily by \eqref{by_above}. 
It remains to verify the lower bound for \(n\,[\psi]\). 
For this we use \eqref{by_below}, noticing that \(\psi\) is 
monotone and that \(|v|^2+|v_*|^2 
\ge |v'|^2+|v'_*|^2\). We then obtain:
\begin{equation*}
\begin{split}
n\,[\psi] 
& \ge \psi(|v'|^2+|v'_*|^2) - \psi(|v'|^2) - \psi(|v'_*|^2) 
\\
& \ge b\,|v'|^2 \, |v'_*|^2 \,\psi''(|v'|^2 + |v'_*|^2)
\\
& = b\,\zeta(v',v'_*)\, (|v'|^2 + |v'_*|^2)^2 \, 
\psi''(|v'|^2 + |v'_*|^2)\,,
\end{split}
\end{equation*}
where 
\[
\zeta(v',v'_*) = 
\frac{|v'|^2}{|v'|^2 + |v'_*|^2}\,
\frac{|v'_*|^2}{|v'|^2 + |v'_*|^2}. 
\]
Further, using \eqref{ine-mu}, we get
\begin{equation*}
\zeta(v',v'_*) = \frac{1}{4}
\Big(1-\frac{4\lambda^2|u|^2|w|^2}{(\lambda^2|u|^2 + |w|^2)^2}\,
\cos^2\mu\Big)
\ge \frac{1}{4}\, (1-\cos^2\mu) = \frac{1}{4}\,\sin^2\!\mu . 
\end{equation*}
Finally, noticing that
\begin{equation*}
|v'|^2+|v'_*|^2 = \frac{\lambda^2|u|^2 + |w|^2}{2}\ge 
\lambda^2\frac{|u|^2 + |w|^2}{2} = \lambda^2 \,(|v|^2+|v_*|^2).
\end{equation*}
we obtain
\[
n\,[\psi] \ge 
\frac{b}{4}\,\lambda^2\,(\eta_2(\lambda^{-2}))^{-1}\,\sin^2\!\mu\,\,
(|v|^2 + |v_*|^2)^2 \, \psi''(|v|^2 + |v_*|^2). 
\]
This completes the proof of the lemma.
\end{proof}

Lemma \ref{povz:conv} gives us the basic formulation of the 
{\it Povzner inequality} for the considered class of test 
functions \(\psi\). In the example \(\psi(x)=x^p\) we have 
\(p\,[\psi] \sim C(|v|^2|v_*|^{2p-2} +|v_*|^2|v|^{2p-2})\)
and \(n[\psi] \sim c(|v|^{2p}+ |v_*|^{2p})\), outside the 
set where \(\kappa(\lambda,\mu)\) is small (which amounts 
to a small set of angles). This implies that that the nonpositive 
term \(-n\,[\psi]\) is dominating, at least when \(|v|>>|v_*|\)
or \(|v_*|>>|v|\), which are the most important regions of 
integration from the point of view of calculation of moments 
(cf. also \cite{De,MiWe}). We can further simplify the 
inequalities and get rid of the dependence on the angular
variables, by integration with respect to \(\sigma
\in S^{\di-1}\). We then obtain the following lemma.  
\begin{lemma}
\label{cor1}
Assume that the function \(\psi\) 
satisfies \eqref{cond:psi-1}--\eqref{cond:psi-4}.
Then 
\[
\begin{split}
\bar{q}\,[\psi]
\le -\, k \, ( |v|^2 + |v_*|^2 )^2\, 
\psi''( |v|^2 &  + |v_*|^2)
 + \,A\, ( |v|^2\,\psi'(|v_*|^2) + |v_*|^2\,\psi'(|v|^2)\,). 
\end{split}
\]
where the constant \(A\) is as in Lemma \ref{povz:conv},
and \(k>0\) is a constant that depends on the function 
\(\psi\) but not on \(\alpha\).
\end{lemma}
\begin{proof} For the proof we notice that 
\(\lambda(\cos\chi)\) is pointwise decreasing as 
\(\alpha\searrow0\) and so, 
\[
\lambda(\cos\chi) \ge \cos\chi, \quad \text{for}
\quad \cos\chi >0 , 
\]
for all \(\alpha>0\). 
We then denote \(\cos\theta = (\nu\cdot\sigma)\), 
\(\;b_0(\cos\theta) = b(u,\sigma)\), and estimate the integral 
\begin{equation}
\label{int:ang}
\is \kappa(\lambda,\mu)\, b_0(\cos\theta)\,d\sigma
\ge \int_{\{\cos\chi>\eps_0,\;\;\sin\mu>\eps_1,\;\;
1-\cos\theta>\eps_2\}} 
\kappa(\lambda,\mu) \,b_0(\cos\theta)\,d\sigma,
\end{equation}
setting \(\eps_0\), \(\eps_1\) and \(\eps_2\) small enough. 
The integrand on the right-hand side of \eqref{int:ang} is 
bounded below by a constant, and so is the area of the 
domain of integration. (The verification of the last statement 
for the condition \(\sin \mu>\eps_1\) is somewhat tedious 
and  is achieved by changing the variables of integration
from  \(\omega\) to \(\sigma\): we omit the technical 
details). We therefore find that the integral \eqref{int:ang}
is bounded below by a constant \(k>0\), independent on 
\(\alpha\). The rest of the claim is easy to verify.
\end{proof}

Finally, we present estimates for the integral 
expression \eqref{qpsi_bar} multiplied by the relative speed,
in the cases when \(\psi(x)\) is given by one of 
the functions \eqref{pow}. 
\begin{lemma}
\label{povz:pow}
Take \(p>1\) and \(\psi(x)=x^p\). Then
\[
|u|\,\bar{q}\,[\psi](v,v_*) \le - k_p (|v|^{2p+1}+|v_*|^{2p+1}) 
+ A_p (|v||v_*|^{2p} + |v|^{2p}|v_*|). 
\]
Also, take \(\psi(x)=(1+x)^p-1\), then
\[
|u|\,\bar{q}\,[\psi](v,v_*) \le 
- k_p (\jap{v}^{2p+1}+\jap{v_*}^{2p+1}) 
+ A_p (\jap{v}\jap{v_*}^{2p} + \jap{v}^{2p}\jap{v_*}). 
\]
Here the constants \(k_p\) and \(A_p\) are independent 
on the restitution coefficient \(\alpha\). 
\end{lemma}
\begin{proof}
We use Lemma \ref{cor1} and the inequalities
\[
||v|-|v_*||\le |u| = |v-v_*| \le |v|+|v_*|.  
\] 
Then in the case \(\psi(x)=x^p\) the bounds have the form
\[
-\,p(p-1)k_p\,|u|\,(|v|^2+|v_*|^2)^p 
+ pA_p\,|u|\,(|v|^2|v_*|^{2p-2} + |v|^{2p-2}|v_*|^{2}) \, .
\]
The terms appearing with the negative sign
are estimated using the inequality
\[
|u|\,(|v|^2+|v_*|^2)^p \ge \frac{1}{2}(|v|^{2p+1}+|v|^{2p+1}) 
- \frac{1}{2}(|v||v_*|^{2p}+|v|^{2p}|v_*|)\, . 
\]
For the remaining terms we have
\[
|u|\,(|v|^2|v_*|^{2p-2} + |v|^{2p-2}|v_*|^{2}) 
\le C_p(|v||v_*|^{2p} + |v|^{2p}|v_*|)\, ,  
\]
which completes the proof of the first part of the lemma. 
The case \(\psi(x)=(1+x)^p-1\) can be treated by arguing 
along the same lines, by using the inequalities
\[
(|v|^2+|v_*|^2)^2 \ge \frac{1}{2}(1+|v|^2+|v_*|^2)^2 - 1
\]
and
\(
|v|\ge (1+|v|^2)^{1/2}-1. 
\)
\end{proof}

\subsection{Estimates for higher-order moments} The Povzner-type 
inequalities of Lemma \ref{povz:pow} allow us to study the topics 
of propagation and appearance of moments.  We find that results
known for the classical Boltzmann equation with ``hard-forces''
interactions \cite{El,De} transfer to present case. We introduce 
the notation
\[
Y_s(t) = \ir f \jap{v}^s\,dv, 
\]
and denote by \(\bar{Y}_s\) the corresponding steady moment.  
\begin{lemma}
\label{moment-ineq}
Let \(f\) be a sufficiently regular and rapidly decaying solution 
of \eqref{be:tdep}. Then, the following differential inequality 
holds:
\begin{equation}
\label{deq:mom}
\frac{d}{dt}Y_s + 2k_s Y_{s+1} \le K_s ( Y_{s} + Y_{s-2} )
\end{equation}
where \(K_s\) and \(k_s\) are positive constants. Further,
\[
\sup\limits_{t>0} Y_s(t) 
\le Y^*_s = \max\big\{Y_s(0),\big(K_s/k_s\big)^{s}\big\},
\] 
and for every \(\tau>0\)
\begin{equation}
\label{mom:int}
\int_0^\tau Y_{s+1}(t)\,dt \le \frac{K_s\tau+1/2}{k_s}Y^*_s.
\end{equation}
Finally, for the steady 
equation \eqref{be:stat} we obtain the apriori estimate
\[
\bar{Y}_{s+1} \le 
\frac{K_s}{2k_s} (\bar{Y}_{s} + \bar{Y}_{s-2}).
\] 
\end{lemma}

\begin{proof}[Proof of Lemma \ref{moment-ineq}.]
Using the weak form of equation(\ref{be:tdep}) with 
\(\psi(v)=\jap{v}^s\) we find
\begin{equation}
\label{mom2}
\frac{d}{dt} \ir f \, \jap{v}^s \, dv - {\ir} \Delta f\,\jap{v}^s\,dv 
= {\ir} Q(f,f)\,\jap{v}^s\,dv. 
\end{equation}
Estimating the moments of the collision integral according to 
Lemma \ref{povz:pow} we get
\begin{equation*}
{\ir} Q(f,f)\,\jap{v}^s\,dv 
\le - 2 k_s Y_{s+1} + 2 A_s Y_1 Y_s.
\end{equation*}
The moments of the Laplacian term are computed as follows: 
\begin{equation}\label{lapl_mom} 
\begin{split}
\ir f \Delta\jap{v}^s \,dv 
= \ir 
& f (\,  (s(s-2)+s\di) \jap{v}^{s-2} - s(s-2)\jap{v}^{s-4}\, ) \, dv
\\
&=(s(s-2)+s\di) \,Y_{s-2} - s(s-2)\, Y_{s-4}. 
\end{split}
\end{equation}
Combining \eqref{mom2} and \eqref{lapl_mom} and neglecting the 
non-positive \(Y_{s-4}\) term, we obtain inequality \eqref{deq:mom} with 
\(K_s=\max\{2A_s,s(s-2)+s\di\}\). 

\noindent
To obtain a uniform bound for \(Y_s(t)\), we use Jensen's 
inequality to write
\[
Y_{s+1} \ge (Y_{s})^{(s+1)/s}. 
\]
Then we find, estimating the right-hand side of \eqref{deq:mom} 
by \(2K_sY_s\),
\[
\frac{d}{dt}Y_s \le - 2k_s (Y_s)^{(s+1)/s} + 2K_s Y_{s}. 
\]
Thus, \(Y'_s(t)<0\) if \(Y_{s} > (K_s/k_s)^s\), and so, the 
upper bound for \(\sup_{t>0}Y_s(t)\) must hold. 

\noindent
Further, integrating in time we obtain
\[
2k_s \int_0^\tau Y_{s+1} \le 2K_s\tau Y^*_s - Y(s) + Y(0) 
\le (2K_s\tau + 1) Y^*_s, 
\]
which proves \eqref{mom:int}. 

\noindent
Finally, the last inequality is obtained by the same arguments 
as \eqref{deq:mom} applied to the steady equation. 
\end{proof}

Based on the Lemma just proven we can make the following
conclusions about the behavior of the moments of the solutions. 
First, if a moment \(Y_s\) is finite initially, 
it {\it propagates}, that is, it remains bounded for the whole 
time-evolution. Further, the integral condition on 
\(Y_{s+1}\) implies the {\it appearance} of moments
of order \(s+1\): these moments become finite after arbitrarily 
short time, even if they are initially infinite (cf. \cite{De}). 
Indeed, suppose that \(Y_{s+1}(0)=+\infty\), then for every 
\(\tau>0\) there is a \(t_0<\tau\) such 
that \(Y_{s+1}(t_0)<+\infty\). Then, applying the Lemma 
to \(Y_{s+1}\), starting with \(t=t_0\), we obtain that 
{\it for every} \(t_0>0\),
\[
\sup_{t>t_0}Y_{s+1}(t) \le C_{t_0,s}, 
\]
which implies the above statement. 
The last part of the Lemma implies
an important statement concerning the moments of the steady 
solution: on the formal level, {\it every} solution that 
has a finite moment of order 
\(s>2\) has finite moments of {\it all positive orders}. In fact, 
in view of the \(L^1_3(\R^\di)\) estimate of the previous 
section, this implies that every solution 
with finite mass should have this property. 

\section{\(L^p\) bounds and apriori regularity}
\label{sec:reg}

In this section we study the apriori regularity of solutions to 
\eqref{be:tdep} and \eqref{be:stat}. The presence of the diffusion 
term in the equation makes it plausible that solutions to 
the steady equation should be smooth, and those for the 
time-dependent equation should gain smoothness after arbitrarily 
short time. However, to realize this idea we {need to} make use 
of the particular structure of the collision term. As we will 
see below, the moment bounds of the previous section will also 
be of crucial importance. We start by establishing the bounds 
for the collision operator in the spaces \(L^p\) with a 
polynomial weight, extending the results well-known in the 
case of the classical Boltzmann equation, {and first
derived by Gustafsson}~\cite{Gu}.  
{Below, we shall establish these bounds by adapting 
the simple strategy that was suggested 
in~\cite[Chapter~2, Section~3.3]{V} and later developed
in~\cite{MV} to establish improved $L^p$ bounds in the 
elastic case.}

\subsection{\(L^p\) bounds for the collision operator}

We will use the following weighted \(L^p\) spaces: 
\begin{equation*}
L^p_k(\R^\di)=\{f\:|\:f\jap{v}^{k}\in{L^p}(\R^\di)\},
\end{equation*}
where \(\jap{v}=(1+|v|^2)^{1/2}\). The necessity to 
introduce a weight comes from the presence of the factor
\(|u|\) in the hard sphere collision term 
\eqref{co:weak-sym}. The collision operator is generally 
unbounded on \(L^p\): in order to control its norm we will 
invoke the \(L^p_k\) norms with higher powers of 
\(\jap{v}\).  The precise formulation of this statement 
is given in next lemma.  
\begin{lemma}
\label{lem:2} 
For every \(1\le p \le \infty\) and every \(k \ge 0\),
\begin{equation*}
\| Q (g,f)\|_{L^p_k} \le C \, 
\big(  \|g\|_{L_{k+1}^p} \|f\|_{L_{k+1}^1}
+ \|g\|_{L_{k+1}^1} \|f\|_{L_{k+1}^p} \big)\,,
\end{equation*}
where \(C\) is a constant depending on \(p\), \(k\) and \(N\) only. 
\end{lemma}
\begin{proof}
We fix an exponent \(1\leq{p}\leq{\infty}\). It is easy to estimate 
the ``loss'' part \(Q^-(g,f) = (g*|v|)f\), using the inequality
\begin{equation*}
|\,g*|v|\,|  \le  \|g\|_{L_1^1} \jap{v},
\end{equation*}
from which it follows
\begin{equation}\label{eq4.1}
\|Q^-(g,f)\|_{L_{k}^p} \le  \|g\|_{L_{1}^1}\|f\|_{L_{k+1}^p}.
\end{equation}

{We now turn to estimate the $Q^+$ term: starting from the weak form 
\eqref{co:weak-sym}, we find}
\begin{equation}
\label{weak_qplus}
\begin{split}
\| Q^+(g,f)\jap{v}^k\|_{L^p} = 
\sup\limits_{\| \psi\|_{L^{p'}} =1} 
\ir{Q^+}(g,f)\,\psi\,\jap{v}^k\,dv. 
\\
=\sup\limits_{\| \psi\|_{L^{p'}} =1}\ir\ir
f \,g_* |u| \is\psi'\,\jap{v'}^k\,b(u,\sigma) \,d\sigma\,dv\,dv_*.
\end{split}
\end{equation}
By using the inequalities
\(
|u| \le \jap{v} + \jap{v_*}
\)  
and 
\(
\jap{v'}^k \le (\jap{v} + \jap{v_*})^k 
\)
the integral \eqref{weak_qplus} is bounded as
\begin{equation}\label{qplus.bd}
\ir\ir
f \,g_* (\jap{v}+\jap{v_*})^{k+1} 
\is \psi'\,b(u,\sigma)\, d\sigma\,
dv\,dv_*.
\end{equation}
We now see that the problem comes down to estimating the integral 
\[
\Sp[\psi](v,v_*)  = \is\psi'\,b(u,\sigma)\,d\sigma
\]
in either \(L^\infty(\R^\di_v,L^{p'}(\R^\di_{v_*}))\) 
or \(L^\infty(\R^\di_{v_*},L^{p'}(\R^\di_{v}))\). In fact, 
we split \(\Sp[\psi]\) into two parts 
\(\Sp_+[\psi]\) and \(\Sp_-[\psi]\) 
and prove the bounds for each of the parts in the 
respective spaces. We set
\begin{equation*}
\Sp_\pm[\psi] (v,v_*)  = \int_{\{\pm u\cdot\sigma>0\}}\psi'\,
b(u,\sigma)\, d\sigma
\end{equation*} 
and establish the bounds for \(\Sp_+\) and \(\Sp_-\) in 
the following proposition. 
\begin{proposition}
\label{boundtestfcn} The operators
\begin{equation*}
\begin{split}
&\Sp_+ : L^{q}(\R^\di) \to L^\infty(\R^\di_v,L^{q}(\R^\di_{v_*})),\\
& \Sp_- : L^{q}(\R^\di) \to 
L^\infty(\R^\di_{v_*},L^{q}(\R^\di_{v})),
\end{split}
\end{equation*}
are bounded for every \(1 \le q \le \infty \).
\end{proposition}
\begin{proof}
We prove the \(L^q\) bounds by interpolation between \(L^\infty\)
and \(L^1\). The \(L^\infty\) estimates are clear due to the 
boundedness of the domain of integration. 
To check the \(L^1\) bounds we assume without loss of generality 
that \(\psi \ge 0\) and calculate the \(L^1\) norms 
as follows:   
\begin{equation*}
\begin{split} 
\|\Sp_-[\psi](v,v_*)\|_{L^1(\R^\di_{v_*})} 
= & \ir\int_{\{u\cdot\sigma<0\}} 
\psi\,\Big(v+\frac{\beta}{2}\,\big(-u
+|u|\sigma\big)\Big) b(u,\sigma)\,d\sigma\,du
\\
= & \ir \psi(v+z) \is  \!
\frac{ \, b\big(u(z,\sigma),\sigma\big)\,
\chi_{\{(\sigma\cdot u(z,\sigma))<0\}}}
{|J_-(u(z,\sigma),\sigma)|} \,d\sigma \,dz .
\end{split}
\end{equation*}
Here, 
\(
z=v'-v=\frac{\beta}2 
\left( -{u} + |u|\sigma\right),
\)
and \(J_-(u,\sigma)\) is the Jacobian of the 
transformation \(u\mapsto{z}\) (for fixed \(\sigma\)):
\[
J_-(u,\sigma)= \Big( \frac{\beta}{2} \Big)^\di 
\Big(-1+\frac{(u\cdot\sigma)}{|u|}\Big) ,
\]
The condition \({u\cdot\sigma}<0\) ensures that \(|J_-|\) 
is bounded below by \(\left( \frac{\beta}2 \right)^{\di}\), and then,
\begin{equation*}
\|\Sp_-[\psi](v,v_*)\|_{L^1(\R^\di_{v_*})} 
\le \Big( \frac{\beta}2 \Big)^{-\di}\is b(u,\sigma)\,d\sigma\, \|\psi\|_{L^1}
= \Big( \frac{\beta}2 \Big)^{-\di}\|\psi\|_{L^1}\,,
\end{equation*}
for every \(v\in\R^\di\).

Similarly, for the \(\Sp_+\) term we have 
\begin{equation*}
\begin{split} 
\|\Sp_+[\psi](v,v_*)\|_{L^1(\R^\di_{v})} 
= & \ir \int_{\{u\cdot\sigma>0\}} 
\psi\,\Big(v_* + \frac{1}{2}\,\big((2-\beta)\,{u} 
+ \beta\,|u|\sigma\big)\Big)\,b(u,\sigma) \, d\sigma\, du 
\\
= & \ir \psi(v+z) \is  \!
\frac{ \, b\big(u(z,\sigma),\sigma\big)\,
\chi_{\{(\sigma\cdot u(z,\sigma))>0\}}}
{|J_+(u(z,\sigma),\sigma)|} \,d\sigma \,dz \,,
\end{split}
\end{equation*}
where now 
\(
z=v'-v_*=\frac{1}{2}\big((2-\beta)){u}+\beta\,|u|\sigma\big),
\)
and 
\[
J_+(u,\sigma)= \Big( \frac{2-\beta}2 \Big)^{{\di}} 
\Big(1+\frac{\beta}{2-\beta}\frac{(u\cdot\sigma)}{|u|}\Big).
\]
Then, since \((u\cdot\sigma)>0\), we can argue similarly 
to the previous case to obtain 
\begin{equation*}
\|\Sp_+[\psi](v,v_*)\|_{L^1(\R^\di_{v})} 
\le \Big( \frac{2-\beta}2 \Big)^{-\di}\|\psi\|_{L^1} ,
\end{equation*}
uniformly in \(v_*\in\R^\di\). The statement of the proposition  
now follows by the Marcinkiewicz interpolation theorem. 
\end{proof}
\noindent
{\it End of proof of Lemma \ref{lem:2}. }
Combining the bound \eqref{qplus.bd} with the ones proven 
in Proposition \ref{boundtestfcn} we find
\begin{equation*}
\begin{split} 
&\ir Q(g,f) \, \psi\, dv 
\\
&\le \ir\ir
f \,g_*\,(\jap{v}+\jap{v_*})^{k+1}\, 
\big(\Sp_+[\psi](v,v_*) + \Sp_-[\psi](v,v_*) \big)\,{dv}\,{dv_*}
\\
&\le C_k \ir g_* \ir
f\, (\jap{v}^{k+1}+\jap{v_*}^{k+1})\, 
\Sp_+[\psi](v,v_*)\,{dv}\,{dv_*}
\\
&+ C_k \ir
f \ir g_*\,(\jap{v}^{k+1}+\jap{v_*}^{k+1})\, 
\Sp_-[\psi](v,v_*)\, {dv_*}\,{dv}
\\
&\le C \,\big(\, \|g\|_{L^1} \|f\|_{L^p_{k+1}} 
+ \|g\|_{L^1_{k+1}} \|f\|_{L^p} + \|f\|_{L^1} \|g\|_{L^p_{k+1}} 
+ \|f\|_{L^1_{k+1}} \|g\|_{L^p}\big).  
\end{split}
\end{equation*}
since \(\|\psi\|_{L^{p'}}=1\). From this the conclusion 
of the lemma follows easily. 
\end{proof}

\subsection{\(H^1\) regularity: Steady state equation}

We start by establishing apriori estimates for solutions 
to the steady equation \eqref{be:stat}, for which the analysis
is performed in a rather more direct way than for the 
time-dependent problem. We first show the bounds in 
the Sobolev spaces with the weight \(\jap{v}^k = 
(1+|v|^2)^{k/2}\): 
\[
H^1_k(\R^\di) = \{f\in L^2_k(\R^\di)\;|\;\nabla f \in L^2_k(\R^\di)\}.
\]
The main tools are the coercivity of the diffusion part, the 
estimates of the collision operator in \(L^p\), and the
interpolation inequalities for \(L^p\) spaces. The constants 
in the estimates are expressed in terms of the \(L^1\) moments.
{ In all this section, we shall assume for simplicity that
$N\geq 3$, but there is no difficulty to adapt the proofs
to cover the case $N=2$ as well.} 
We begin with an estimate for the gradient in \(L^2\). 
\begin{lemma}
\label{thm:steady1}
Assume that the function \(f\in H^1(\R^\di)\cap L^1_\rbeta(\R^\di)\),
where \(\rbeta=\frac{\di+2}{4}\), is a solution of \eqref{be:stat}. 
Then
\begin{equation*}
\|\nabla f\|_{L^2} \le C A \,B^\rbeta,
\end{equation*} 
where 
\[
A=\|f\|_{L^1_\rbeta}, \quad B=\|f\|_{L^1_1}, 
\]
and \(C\) is a constant depending on the dimension. 
\end{lemma}
\begin{proof}
Multiplying  equation (\ref{be:stat}) by \(f\), integrating 
and applying H\"older's inequality  yields 
\begin{equation}\label{energy:weak}
\ir |\nabla f|^2 \,dv = \ir Q(f,f)f \,dv \le C\|f\|_{L^p}
\|Q(f,f)\|_{L^{p'}},
\end{equation}
for all \(1 \le p \le \infty\). 
We choose \(p=\twos=2\di/(\di-2)\), where \(2^*\) is the 
critical Sobolev exponent,  and apply the 
Sobolev's embedding inequality 
\begin{equation}\label{H1.2}
\|f\|_{L^{\twos}} \le C\|\nabla f\|_{L^2}
\end{equation}
(Note: for \(\di=3\),  \(\twos=6\) 
and \((\twos)'=6/5\).) Then, by Lemma~\ref{lem:2}, 
\begin{equation}\label{H1.3}
\|Q(f,f)\|_{L^{(\twos)'}} \le C \|f\|_{L_1^{(\twos)'}}  \|f\|_{L_1^1}.
\end{equation}
 
To estimate \(\|f\|_{L_1^{(\twos)'}}\) we  use the following 
interpolation inequality for weighted \(L^p\) norms (\(\varphi\) 
is any weight function), which can be easily verified 
using H\"older's inequality: 
\begin{equation}\label{interp}
\|f\varphi^k\|_{L^q} \le \|f\varphi^{k_1}\|_{L^{q_1}}^\nu 
\|f\varphi^{k_2}\|_{L^{q_2}}^{1-\nu},
\end{equation}
where 
\begin{equation*}
\frac{\nu}{q_1} + \frac{1-\nu}{q_2} = \frac1q
{\quad\rm{}and\quad } k_1\nu + k_2 (1-\nu) = k.
\end{equation*}
Now, interpolating the norm in \(L^q_1\) for $q={(\twos)'}$ between  
$q_1 = \twos$ and  $q_2 =1$, we get  
\begin{equation}\label{f.twos.prime}
\|f\|_{L_1^{(\twos)'}} \le \|f\|_{L^{\twos}}^\nu 
\|f\|_{L^{1}_\rbeta}^{1-\nu}
\end{equation}
where \(\nu\) and \(\rbeta\) are determined from the following equations:
\begin{equation*}
\frac{\nu}{\twos} + \frac{1-\nu}1 = \frac1{(\twos)'}\quad \text{and}\quad 
\rbeta\left( 1-\nu\right) =1,
\end{equation*}
so that 
\begin{equation}\label{nu.and.beta}
\nu = \frac{\di-2}{\di+2}\quad \text{and}\quad 
\rbeta = \frac{1}{1-\nu} = \frac{\di+2}4.
\end{equation}
Combining estimates  \eqref{energy:weak}--\eqref{f.twos.prime}
we obtain the inequality
\begin{equation}\label{H1.4}
\begin{split} 
\|\nabla f\|_{L^2}^2 
\le C \|f\|_{L_1^1} \|f\|_{L_{\rbeta}^1}^{1-\nu}\|f\|_{L^{2*}}^{1+\nu} 
\le C B A^{1-\nu} \|\nabla f\|_{L^2}^{1+\nu},
\end{split}
\end{equation}
from which the conclusion of the lemma follows. 
\end{proof}

The result of the Lemma implies a bound for the solutions in 
the space \(H^1(\R^\di)\). Indeed, by the Sobolev embedding,
\[
\|f\|_{L^\twos}\le C A \,B^\rbeta. 
\]  
Interpolating between \(L^1\) and \(L^\twos\) using inequality
\eqref{interp} we get a bound for the \(L^2\) norm, which then 
implies a bound in \(H^1\). Since the constants in the estimates 
depend on the  \(L^1_k\) norms only, and the latter are controlled 
by the moments bounds, we gain an apriori control of the \(H^1\)
norm by means of the mass and the energy only.  We next 
see that the derivatives of the solutions have an appropriate 
decay, so even \(L^2_k\) norms for all \(k\ge0\) are bounded.   

\begin{lemma}
\label{le:l2mom}
Let \(f\) be a solution of equation \eqref{be:stat} and assume 
that \(f\in H^1_k(\R^\di)\,\cap\, L^1_{(k+1)\rbeta}(\R^\di)\), 
where \(k\ge0\) and \(\rbeta=\frac{\di+2}{4}\). Then 
\[
\|\nabla (f\jap{v}^k)\|_{L^2} \le C\big(A_1 \,A_2^\rbeta
+k^{2\rbeta} A_3\big), 
\]
where 
\[
A_1=\|f\|_{L^1_{(k+1)\rbeta}}, \quad A_2 = \|f\|_{L^1_{k+1}}, 
\quad A_3=\|f\|_{L^1_{k-2/\rbeta}},
\] 
and \(C\) is a constant 
depending on the dimension \(\di\).
\end{lemma}
\begin{proof}
Integrating equation \eqref{be:stat} against 
\(f\jap{v}^{2k}\) we obtain 
\begin{equation}\label{eq:gradmom}
\ir\nabla f\cdot \nabla (f\jap{v}^{2k})\,dv 
= \ir Q(f,f) f\jap{v}^{2k}\,dv.
\end{equation}
Using estimates from the previous lemma,
the right-hand side can be bounded above as follows:
\begin{equation}\label{rhs1}
\begin{split}
& \|Q(f,f)\jap{v}^k\|_{L^{(\twos)'}} 
\|f\jap{v}^k\|_{L^{2*}}\\
\myskip 
& \le C\|f\jap{v}^{k+1}\|_{L^{(\twos)'}} 
\|f\jap{v}^{k+1}\|_{L^1} 
\|f\jap{v}^k\|_{L^{\twos}}. 
\end{split}
\end{equation}
Interpolating as in \eqref{interp} we find
\begin{equation}\label{rhs2}
\begin{split}
\|f\jap{v}^{k+1}\|_{L^{(\twos)'}} 
\le \|f\|_{L^{2*}}^\nu 
\|f\|_{L^1_{(k+1)\rbeta}}^{1-\nu},
\end{split}
\end{equation}
where \(\nu\) and \(\rbeta\) are as defined in \eqref{nu.and.beta}.
Therefore, combining  \eqref{rhs1} with \eqref{rhs2} 
we bound the right hand side of  \eqref{eq:gradmom} by
\begin{equation}\label{rhs3}
\begin{split}
&C\big(\, \|f\|_{L^1_{(k+1)\rbeta}}^{1-\nu}
\|f\|_{L^1_{k+1}}\big)\|f\jap{v}^k\|_{L^{2*}}^{1+\nu}\\ 
&\le C A_1^{1-\nu} A_2 \|\nabla (f\jap{v}^k)\|_{L^2}^{1+\nu}.
\end{split}
\end{equation}
The integral on the left-hand side 
of \eqref{eq:gradmom} is estimated as follows:
\begin{equation}\label{lhs:est}
\begin{split}
\ir \nabla f \cdot \nabla (f\jap{v}^{2k}) \, dv 
=\ir |\nabla & (f\jap{v}^{k}) |^2 \, dv - 
\ir f^2 |\nabla\jap{v}^k|^2\,dv
\\
&\ge \| \nabla (f\jap{v}^{k}) \|^2_{L^2} - k^2\|f\|^2_{L^2_{k-1}}. 
\end{split}
\end{equation}
Further, interpolating the \(L^2_{k-1}\) norm 
between \(L^{\twos}\) and \(L^1\) we get 
\begin{equation*}
\begin{split}
\|f\jap{v}^{k-1}\|_{L^2}  
\le  C\|f\jap{v}^k\|_{L^{2*}}^\lambda 
\|f\jap{v}^{k_2}\|_{L^1}^{1-\lambda} 
\le C \|\nabla (f \jap{v}^k)\|_{L^2}^\lambda
\|f\|_{L^1_{k_2}}^{1-\lambda}, 
\end{split}
\end{equation*}
where
\begin{equation}\label{interp:lambda}
\frac{\lambda}{\twos} + \frac{1-\lambda}1 
= \frac12,\quad\text{so that}\quad \lambda = \frac{\di}{\di+2} 
= \frac{1+\nu}{2},
\end{equation}
and
\[
k-1 = \lambda k + (1-\lambda) k_2, 
\quad\text{so that}\quad 
k_2 = k - \frac{1}{1-\lambda} = k - 2/\rbeta .
\]
Gathering the above inequalities and noticing that \(2\lambda=1+\nu\) 
we obtain:
\[
\|\nabla (f \jap{v}^k)\|_{L^2}^2
\le
C \big(A_1^{1-\nu} A_2   
+ k^2 A_3^{1-\nu}\big)\|\nabla (f \jap{v}^k)\|_{L^2}^{1+\nu}.
\]
Dividing by the norm of the gradient to the power \(1+\nu\) we get
\[
\|\nabla (f \jap{v}^k)\|_{L^2}\le \big( C A_1^{1-\nu} A_2 + 
k^2 A_3^{1-\nu}\big)^\frac{1}{1-\nu}.
\]
Noticing that \(\frac{1}{1-\nu}=\rbeta\) and using the inequality
\((x+y)^r\le C_r(x^r+ y^r)\)
we arrive at the 
conclusion of the lemma.  
\end{proof}

Using the Lemma just proven we find bounds for solutions \(f\) 
in \(H^1_k\) for every \(k \ge 0 \). Indeed, using the inequality
\[
|(\nabla f)\jap{v}^k|^2 
\le C(\,|\nabla (f \jap{v}^k) |^2 + |f\nabla\jap{v}^k|^2)
\]
and interpolating in the second term between 
\(L^{\twos}\) and \(L^1_{k-2/\rbeta}\) we get 
\[
\|\nabla f\|^2_{L^2_k} 
\le C \big( \|\nabla(f \jap{v}^k)\|^2_{L^2} 
+ \|\nabla(f\jap{v}^k)\|_{L^{2}}^{1+\nu} 
\|f\|_{L^1_{k-2/\rbeta}}^{1-\nu} \big),
\]
from which an estimate in terms of the \(L^1\) moments 
follows.  Further, by interpolation inequality \eqref{interp}, 
\[
\|f\|_{L^2_k} \le \|f\|_{L^\twos}^\lambda 
\|f\|^{1-\lambda}_{L^1_{k/(1-\lambda)}}, 
\]
and so, in view of our earlier remarks, the norm in 
\(L^2_k\) is also estimated in terms of \(L^1\) moments 
only. Summarizing the results obtained so far, the solutions 
are controlled apriori in \(H^1_k(\R^\di)\) for any \(k\ge0\) 
in terms of mass and kinetic energy only.   

\subsection{Schwartz class regularity: Steady problem}
Our next aim is now to establish a priori bounds for solutions to 
\eqref{be:stat} in the spaces 
\[
H^n_k(\R^\di)=\{f\in L^2_k(\R^\di)\; |\; 
\nabla^m{f} \in L^2_k(\R^\di), \; 1 \le m \le n\},
\] 
for all \( 1 \le n < \infty \) and all \( 0 \le k < \infty \). 
We use induction on \(n\), differentiating 
the equation in \(v\) in each step. The base of the induction 
is given by Lemma \ref{le:l2mom}. 
We recall the following rule for differentiating the collision
integral. 
\begin{proposition}
\label{diff_Q}
Let \(f\) and \(g\) be smooth, rapidly decaying 
functions of \(v\). Then
\[
\nabla{Q(g,f)}=Q(\nabla{g},f)+Q({g},\nabla{f}). 
\]
\end{proposition}
\begin{proof} We use the splitting into the ``gain'' and 
``loss'' terms, \(Q(g,f)=Q^{+}(g,f)-Q^-(g,f)\). 
Since \(Q^-(g,f)=f(g*|v|)\), the differentiation rule for 
the ``loss'' term is obvious. To prove the proposition for 
the ``gain'' term \(Q^+(g,f)\) we represent it as follows,
using \eqref{col:pre}:  
\begin{equation*}
{Q^+}(g,f)
=\ir\is 
f\Big(v+\frac{\lprime{u}-u}{2}\Big)
g\Big(v-\frac{\lprime{u}+u}{2}\Big)
\,\frac{1}{\alpha^2}\,|u|\,b(u,\sigma)\,d\sigma\,du .
\end{equation*}
Since \(\lprime{u}\) is a function of \(u\) and \(\sigma\)
only, the statement follows by differentiation under the 
integral sign. 
\end{proof}

\begin{remark}
The above statement is in fact a corollary of the following
abstract statement which can be proven very easily: Let $Q$ 
be a bilinear operator commuting with translations, continuously 
differentiable; then $\nabla Q(g,f)= Q(\nabla f, g) + 
Q(f,\nabla g)$. Thus, the differentiation formula 
of Proposition \ref{diff_Q} can be seen as a consequence 
of the translation invariance of \(Q\). 
\end{remark}

As a direct corollary of Proposition \ref{diff_Q}, 
higher-order derivatives of \(Q\) can be 
calculated using the following Leibniz formula: 
\begin{equation*}
\partial^j {Q}(g,f)=\sum\limits_{0\le l \le j}
{\lchoose{j}{l}}
\,Q(\partial^{j-l}g,\partial^{l}{f}),
\end{equation*}
where \(j\) and \(l\) are multi-indices \(j=(j_1\dots j_\di)\), and 
\(l=(l_1\dots l_\di)\); 
\[
\partial^j = \partial^{j_1}_{v_1} \dots \partial^{j_\di}_{v_\di}, 
\]
and \(\lchoose{j}{l}\) are the multinomial coefficients. 
Thus, for every multi-index \(j\), by formal differentiation 
of \eqref{be:stat} we obtain the following equations 
for higher-order derivatives: 
\begin{equation}
\label{be:diff}
-\Delta \, \partial^j {f}=\sum\limits_{0\le l \le j}
{\lchoose{j}{l}}\,Q(\partial^{j-l}g,\partial^{l}{f}).
\end{equation}
By applying the methods developed in Lemmas \ref{thm:steady1} 
and \ref{le:l2mom} to equation \eqref{be:diff} we 
arrive at the following result. 
\begin{lemma}\label{le:l2diff}
Let \(f\) be a solution to \eqref{be:stat}, 
such that \(f\in H^{n+1}_{k+\mu}(\R^\di)\), with
\(n \ge 0\), \(k \ge 0\) and \(\mu>1+\frac{N}{2}\). 
Then 
\begin{equation*}
\begin{split}
\|\nabla^{n+1} f\|_{L^2_k} \le 
C \, (1+ k + \|f\|_{H^{n-1}_{k+\mu}})\,
( 1 + \|\nabla^{n} f\|_{L^2_{k+\mu}} ) ,  
\end{split}
\end{equation*}
where \(C\) is a constant depending on \(n\) and  
\(\di\) only. 
\end{lemma}
\begin{proof}
Taking a multi-index \(j\) with \(|j|=n\), multiplying 
equation \eqref{be:diff} by \(\partial^j {f} \, \jap{v}^{2k} \) 
and integrating by parts
we obtain:
\begin{equation}\label{deriv:en}
\begin{split}
\ir \nabla \partial^j f 
& \cdot \nabla (\partial^j f \,\jap{v}^{2k}) \, dv
\\
& = 
\sum\limits_{0\le l \le j}
{\lchoose{j}{l}}\, \ir Q(\partial^{j-l}f,\partial^{l}{f})\, 
\partial^j f \, \jap{v}^{2k} \, dv.
\end{split}
\end{equation}
Similarly to \eqref{lhs:est}, the left-hand side can be written as 
\[
\|\nabla (\partial^j f\, \jap{v}^k  ) \|^2_{L^2} - 
\|\partial^j f \nabla\jap{v}^k\|^2_{L^2}.
\]
Each integral on the right-hand side of \eqref{deriv:en} 
can be bounded above by using Cauchy-Schwartz's inequality 
and Lemma \ref{lem:2} as follows:
\begin{equation*}
\begin{split}
& \ir Q(\partial^{j-l}f,\partial^{l}{f})\, 
\partial^j f \, \jap{v}^{2k} \, dv 
\\
& \le \|Q(\partial^{j-l} f,\partial^l {f})\|_{L^2_k}
\|\partial^j f \|_{L^2_k}
\le C \|\partial^{l} {f}\|_{L^1_{k+1}}
\|\partial^{l-j} f\|_{L^2_{k+1}}
\|\nabla^n f \|^2_{L^2_{k+1}}
\|\partial^j f \|_{L^2_k}.
\end{split}
\end{equation*}
Now, the \(L^1\) norms can be estimated as follows: 
\[
\|\partial^{l} {f}\|_{L^1_{k+1}} 
\le \|\jap{v}^{1-\mu}\|_{L^2}\,
\|\partial^{l} {f}\|_{L^2_{k+\mu}}
\le C\,\|\partial^{l} {f}\|_{L^2_{k+\mu}}
\]
as soon as \(\mu>1+\frac{\di}{2}\).  
Gathering the above estimates we obtain:
\[
\begin{split}
&\|( \partial^j f \, \jap{v}^k ) \|^2_{L^2} 
\\
&\le 
\|\partial^j f \nabla \jap{v}^k\|^2_{L^2} + 
C \, \|\partial^j f\|_{L^2_{k}}
\sum\limits_{0\le l \le j}
{\lchoose{j}{l}}\,  
\|\partial^{l} {f}\|_{L^2_{k+\mu}}\|
\partial^{l-j} f\|_{L^2_{k+1}}.
\\
&\le k^2 \|\partial^j f\|^2_{L^2_{k-1}} 
+ C \|f\|_{L^2_{k+\mu}}  
\|\partial^j f\|^2_{L^2_{k+\mu}} + 
C \|\partial^j f\|_{L^2_{k+\mu}} 
\|f\|^2_{H^{n-1}_{k+\mu}}
\\
&\le k^2 \|\partial^j f\|^2_{L^2_{k+\mu}} + 
C ( 1 + \|f\|^2_{H^{n-1}_{k+\mu}} )  
\|\partial^j f\|^2_{L^2_{k+\mu}} + 
C (1 + \|\partial^j f\|^2_{L^2_{k+\mu}} ) 
\|f\|^2_{H^{n-1}_{k+\mu}}.
\end{split}
\]
Since 
\[
\nabla ( \partial^j f \, \jap{v}^k ) 
= (\nabla \partial^j f )\, \jap{v}^k  + 
 k\,\partial^j f \, |v| \jap{v}^{k-2}, 
\]
we obtain
\[
\|\nabla  \partial^j f \|^2_{L^2_k}
\le C \,(1 +  k^2 + \|f\|^2_{H^{n-1}_{k+\mu}} ) \,
( 1 + \|\partial^j f \|^2_{L^2_{k+\mu}} ) . 
\]
Taking the sum over all \(j\) with \(|j|=n\) implies the estimate 
of the lemma. 
\end{proof}

Lemma~\ref{le:l2diff} gives us a way to estimate higher-order 
derivatives of solutions in terms of lower-order ones. 
Thus, provided we have a solution to \eqref{be:stat} that 
has  all \(H^1_k\) norms bounded in terms of mass and energy
(as we assumed in the previous section), we can derive bounds 
in \(H^2_k\) for every \(k\), and then proceed by induction, 
obtaining bounds in \(H^n_k(\R^\di)\), for all \(n\) and all 
\(k\ge0\). We then obtain
\[
f\in\bigcap\limits_{n \ge 1,\, k\ge 0}H^n_k = \mathcal{S}, 
\]
where \(\mathcal{S}\) is the Schwartz class of rapidly decaying 
smooth functions. Notice that the bounds in each of the 
spaces \(H^n_k(\R^\di)\) can be expressed in terms of mass 
and energy of the solutions.  

\subsection{Regularity for the time-dependent problem}

An analysis of the regularity of the time-dependent solutions
can be performed in the same vein as for the steady problem. 
Using the estimates obtained in the previous section in 
combination with Gronwall lemma will give us results for 
the time-dependent equation \eqref{be:tdep}.  
Our first lemma is an analog of Lemma \ref{thm:steady1}.

\begin{lemma}
\label{reg:tdep-1}
Let \(f\) be a sufficiently regular solution to \eqref{be:tdep} 
with the initial condition \(f(\cdot,0) = f_0\in L^2(\R^\di)\),
such that \(f\) has a moment of order \(\rbeta = \frac{\di+2}{4}\)
bounded uniformly in time. Then
\[
\|f(\cdot,t)\|_{L^2} \le C,\quad 0\le t <\infty,
\]  
and 
\[
\|\nabla f\|_{L^2([0,T]\times\R^\di)} \le C_T, 
\]
for every \(0\le T <\infty\), where the constants \(C\) and \(C_T\)
depend only on \(\di\), \(\|f_0\|_{L^2}\) and 
\(\sup\limits_{t\ge0}\|f(\cdot,t)\|_{L^1_\rbeta}\). 
\end{lemma}
\begin{proof}
Integrating equation \eqref{be:tdep}
against \(f\) we get, arguing similarly to the case 
of the steady problem:
\begin{equation}\label{weak:L2}
\frac{1}{2}\frac{d}{dt}\|f\|^2_{L^2}+\|\nabla{f}\|^2_{L^2}
\le  K(t) \|\nabla f\|_{L^2}^{1+\nu}, 
\end{equation}
where \(K(t)=C\,\|f\|_{L_1^1} \|f\|_{L_{\rbeta}^1}^{1-\nu}\) and 
\(\nu=\frac{\di-2}{\di+2}\) as in \eqref{nu.and.beta}.
By interpolation and Sobolev embedding,
\[
\|f\|_{L^2} \le \|f\|^{\lambda}_{L^{2*}}\|f\|^{1-\lambda}_{L^1}
\le C \|\nabla{f}\|^{\lambda}_{L^{2}}\|f\|^{1-\lambda}_{L^1},
\]
where \(\lambda=\frac{\di}{\di+2}\) as given by \eqref{interp:lambda}. 
Therefore, 
\begin{equation}
\label{reg:eq1}
\|\nabla{f}\|^2_{L^{2}} \ge k \|f\|^{2/\lambda}_{L^2},
\end{equation}
where \(k=C^{-1}\|f\|^{(\lambda-1)/\lambda}_{L^1}\) is a constant. 
Distributing the term \(\|\nabla f\|_{L^2}\) in \eqref{weak:L2}
equally between the left and the right-hand sides and using 
inequality \eqref{reg:eq1} we obtain
\begin{equation}
\label{reg:eq2}
\frac{d}{dt}\|f\|^2_{L^2}+k\|f\|^{2/\lambda}_{L^2}
\le - \|\nabla{f}\|^2_{L^2} + 2K(t) \|\nabla f\|_{L^2}^{1+\nu}
\end{equation}
The function \(X\mapsto-X^2+2K(t)X^{1+\nu}\), appearing 
on the right-hand side of \eqref{reg:eq2} has a global maximum 
\((1+\nu)^{2\rbeta-1}(1-\nu) K(t)^{2\rbeta} = 
C K(t)^{2\rbeta}\), so we obtain 
\begin{equation}
\label{reg:eq3}
\frac{d}{dt}\|f\|^2_{L^2}+k\|f\|^{2/\lambda}_{L^2}
\le C K(t)^{2\rbeta}
\le C \bar{K}^{2\rbeta}. 
\end{equation}
where \(\bar{K}=\sup\limits_{t\ge 0}K(t)\le\sup\limits_{t\ge 0}
\|f\|_{L^1_\rbeta}^2\).  
Applying a Gronwall's lemma argument to \eqref{reg:eq3} 
we then obtain a bound of the  \({L^2}\) norm of \(f\) 
in terms of \(\|f_0\|_{L^2}\) and \(\sup\limits_{t\ge 0}
\|f\|_{L^1_\rbeta}\). Further, integrating \eqref{weak:L2} 
over time, we get
\[
\|\nabla f\|^2_{L^2([0,T]\times\R^\di)} \le C 
+ \bar{K} \int_0^T\|\nabla f\|^{1+\nu}_{L^2}\,dt 
\le C + \bar{K} T^{(1-\nu)/2} \|\nabla f\|^{1+\nu}_{L^2([0,T]\times\R^\di)},
\]
which proves the second claim of the lemma. 
\end{proof}

Similar results can be established about the time-dependence of 
the \(L^2_k\) norms of the solutions. 

\begin{lemma}
\label{reg:tdep-2}
Let \(f\) be a sufficiently regular solution to \eqref{be:tdep},
with initial data \(f_0\in L^2_k(\R^\di)\), where \(k\ge0\), 
and such that \(f\) has a moment of order \(\rbeta(k+1)\), 
where \(\rbeta=\frac{\di+2}{4}\), bounded uniformly in time. 
Then
\[
\|f(\cdot,t)\|_{L^2_k} \le C,\quad 0\le t <\infty,
\]  
and 
\[
\|\nabla f\|_{L^2_k([0,T]\times\R^\di)} \le C_T, 
\]
for every \(0\le T <\infty\), where the constants \(C\) and 
\(C_T\) depend on \(\di\), \(\|f_0\|_{L^2_k}\) and 
\(\sup\limits_{t\ge0}\|f(\cdot,t)\|_{L^1_{\rbeta(k+1)}}\) only. 
\end{lemma}
\begin{proof}
Multiplying the equation by \(f\jap{v}^{2k}\) 
and integrating we obtain:
\begin{equation*}
\begin{split}
\frac{1}{2}\frac{d}{dt}\|f\|^2_{L^2_k} 
&+\ir \nabla f \cdot \nabla (f\jap{v}^{2k}) \, dv 
= \ir Q(f,f) f\jap{v}^{2k}\,dv
\end{split}
\end{equation*}
Following the steps of the proof of Lemma \ref{le:l2mom} and distributing 
the term \(\|\nabla(f\jap{v}^{k})\|^2\) evenly between the left and 
right-hand sides we obtain the following differential inequality:
\begin{equation}\label{dieq}
\begin{split}
&\frac{d}{dt}\|f\jap{v}^k\|^2_{L^2}+\|\nabla(f\jap{v}^{k})\|^2_{L^2}
\le 
-\| \nabla (f\jap{v}^{k}) \|^2_{L^2} 
\\
&
+ 
\left( C A_1(t)^{1-\nu}A_2(t) 
+ k^2 A_3(t)^{1-\nu}\right)
\|\nabla (f\jap{v}^{k})\|_{L^2}^{1+\nu},
\end{split}
\end{equation}
where \(A_1(t)\), \(A_2(t)\), and \(A_3(t)\) are the moments 
defined in Lemma \ref{le:l2mom}. 
The uniform bounds of the moments imply
that the right-hand side of \eqref{dieq} is bounded above by a 
constant. The left-hand side is estimated below as 
\[
\frac{d}{dt}\|f\jap{v}^k\|^2_{L^2}+c\|f\jap{v}^k\|^{2/\lambda}_{L^2}
\]
analogously to \eqref{weak:L2}. Thus, by a Gronwall-type 
argument we obtain that the \(L^2_k\)-norm of \(f\) is bounded 
uniformly in time. Integrating \eqref{dieq} over time we also get 
the second claim of the lemma. 
\end{proof}

Finally, we establish the following analog of Lemma \ref{le:l2diff}
which will allow us to study the regularity of higher-order 
derivatives. 

\begin{lemma}
\label{reg:tdep-3}
Let \(f\) be a solution to \eqref{be:tdep} with initial data 
\(f_0 \in H^n_k(\R^\di)\) where \(k \ge 0\) and \(n \ge 0\), 
such that \(f\) has a moment of 
order \(r^*=\rbeta(2^n (k+\mu)-2\mu+1)\), 
where \(\rbeta=\frac{\di+2}{4}\) 
and \(\mu>\frac{\di+2}{2}\), bounded uniformly in time. 
Then
\[
\|f(\cdot,t)\|_{H^n_k} \le C,\quad 0\le t <\infty,
\]  
and 
\[
\|f\|_{L^2([0,T],H^{n+1}_k(\R^\di))} \le C_T, 
\]
for every \(0\le T <\infty\), where the constants 
\(C\) and \(C_T\) depend on \(\di\), \(\|f_0\|_{H^n_k}\)
and \(\sup\limits_{t\ge0} \|f(\cdot,t)\|_{L^1_{r^*}}\,\) only.   
\end{lemma}
\begin{proof} We will use induction on \(n\). The case \(n=0\) is 
already proven in Lemma \ref{reg:tdep-2}.   Assuming that 
the statement of the lemma holds for \(n-1\), we 
differentiate the equation in \(v\) and argue as in the proof of 
Lemma \ref{le:l2diff}, obtaining the following
inequality:
\[
\frac{1}{2}\frac{d}{dt} \|\nabla^n f\|^2_{L^2_k} + 
\|\nabla^{n+1} f\|^2_{L^2_k} \le C (1+k^2+\|f\|^2_{H^{n-1}_{k+\mu}})
(1+\|\nabla^n f\|^2_{L^2_{k+\mu}})
\] 
We estimate \(\|\nabla^n f\|^2_{L^2_{k+\mu}}\) integrating by parts 
and using Young's inequality (cf. \cite{DeVi}): 
\[
\|\nabla^n f\|^2_{L^2_{k+\mu}} \le \delta \|\nabla^{n+1}f\|^2_{L^2} 
+ C_\delta \|\nabla^{n-1}f\|_{L^2_{2(k+\mu)}}. 
\]
Then, since we assumed \(f\) to be bounded in 
\(H^{n-1}_{2(k+\mu)}\) we get 
\[
\frac{1}{2}\frac{d}{dt} \|\nabla^n f\|^2_{L^2_k} \le
- \|\nabla^{n+1} f\|^2_{L^2_k} + C_1(n,k) \delta \|\nabla^{n+1} f\|^2_{L^2}
+ C_2(n,k,\delta). 
\] 
Choosing \(\delta\) suitably small, we obtain the conclusion 
by Gronwall's lemma. 
\end{proof}

Lemmas \ref{reg:tdep-1} -- \ref{reg:tdep-3} allow 
us to make the following conclusions about the regularity of
solutions to \eqref{be:tdep}. Provided a sufficient number of 
moments is initially available, the \(H^n\) regularity of the 
initial data is preserved with time. Moreover, the established
bounds for the derivatives in \(L^2([0,T]\times\R^\di)\)
imply that after {\it arbitrarily short time} the derivatives 
\(\partial^j f(\cdot,t)\) of {\it any order} are in 
\(L^2(\R^\di)\), and then they propagate in time. Thus,
on the level of apriori estimates we find that the solutions 
become immediately infinitely smooth in \(v\) and decay 
faster than any negative power for \(|v|\) large.   

We can also see that the solutions are infinitely differentiable 
in \(t\). Indeed, in view of the established \(H^n_k\) 
regularity we have 
\(f(\cdot,t)\in \mathcal{S}(\R^\di)\) for \(t>0\), 
and then equation \eqref{be:tdep} implies
\(
\partial_t f(\cdot,t) \in \mathcal{S}(\R^\di), 
\text{ for every } t>0.   
\) 
Differentiating the equation in time and proceeding by 
induction we find also that \(\partial^m_t f(\cdot,t)\in 
\mathcal{S}(\R^\di)\), for every \(m=1,2,\dots\,\), and
for every \(t>0\). The time derivatives also remain 
bounded uniformly in time.   

\section{Existence} 
\label{sec:exist}

We next proceed with a rigorous proof of existence that will also 
justify the formal manipulations performed in the derivation of 
apriori inequalities. 

\begin{theorem}\label{exist:tdep}
For every  \(f_0\ge0\), \(f_0\in L^1_{2}\cap L\log{L}(\R^\di)\)
there exists a nonnegative weak solution
\[
f\in L^\infty( [0,\infty), L^1_{2}(\R^\di)),\quad 
f\log{f} \in L^\infty( [0,\infty), L^1(\R^\di))
\]
to equation \eqref{be:tdep}, with the initial 
condition \(f(\cdot,0)=f_0\). 
Furthermore, if in addition \(f_0\in L^1_{\rbeta}\cap L^2(\R^\di)\),  where 
\(\rbeta=\max\{2,\frac{\di+2}{4}\}\), then 
for every \(t_0 > 0\),   
\[
f\in C^\infty_b ([t_0,\infty), \mathcal{S}(\R^\di)),
\]
where \(C^\infty_b\) denotes the class of functions with 
bounded derivatives of any order, and \(\mathcal{S}\) is 
the Schwartz class of rapidly decaying smooth functions. 
In particular, for \(t>0\), \(f\) is a classical solution 
of \eqref{be:tdep}. 
\end{theorem}
\begin{theorem}\label{exist:stat}
For every \(\rho>0\) there exists a nonnegative 
solution \(f\) to \eqref{be:stat}, 
\[
f\in \mathcal{S}(\R^\di), \quad \text{satisfying}\quad \ir f\,dv=\rho.
\]
Furthermore, every nonnegative solution in
\(
L^1_{\rbeta} \cap L^2(\R^\di)),
\)
where \(\rbeta=\max\{2,\frac{\di+2}{4}\}\) 
is in fact in \(\mathcal{S}(\R^\di)\). 
\end{theorem}

\begin{proof}[Proof of Theorem \ref{exist:tdep}]
We assume that the initial datum \(f_0\) is in 
\(C^\infty(\R^\di)\) and has compact support (we 
will remove this assumption in the end of the 
proof). We also introduce a truncation in the 
collision term by replacing the factor \(|u|\)
in \eqref{co:weak-sym} by 
\begin{equation}
\label{B:trunc}
|u|_{m,M}=m+\min\{|u|,M\}
\end{equation}
and \(m>0\), 
\(M>0\) are truncation parameters. We then 
denote by \(Q_{m,M}(f,f)\) the corresponding 
collision operator. 

The first step of the proof will be to find 
approximating solutions which we define using 
the following truncated problem
\begin{equation}
\begin{split}
\label{fp:trunc}
\partial_t{f}-\Delta_v{f} & = Q_{m,M}(f,f), 
\quad v \in \R^3, \quad t \in [0,T], 
\\
{f}(0,v) & =f_0(v), 
\end{split}
\end{equation}
where \(m\), \(M\) and \(T\) are fixed positive parameters. 
We will denote by \(f\) solutions to \eqref{fp:trunc}, keeping 
in mind that they generally depend on \(m\) and \(M\).

The solutions will be constructed by 
applying a fixed point argument to the following approximation 
scheme:  
\begin{equation}
\begin{split}
\label{fp}
\partial_t{f} - \Delta_v{f} + M{f} 
& = Q_{m,M}(g,g) + M{g}, \quad v \in \R^3, \quad t \in [0,T],
\\
f(v,0) & = f_0(v).
\end{split}
\end{equation}
Here \(g\) is a nonnegative function from 
\(L^\infty([0,T],L^1_2 \cap L^2(\R^\di))\), 
which for every \(t>0\) has unit mass and zero average. 

Denoting by \(h\) the right hand side of equation 
\eqref{fp} we notice that \(h\ge0\), for every \(g\ge0\),  
due to the truncation of the kernel. Indeed, 
\begin{equation}\label{fp:rhs}
h = Q_{m,M}(g,g) + Mg \ge - g \,( g*|v|_{m,M} ) + Mg \ge 0.   
\end{equation}
Further, by analogy with Lemma \ref{lem:2} we can estimate 
\(Q_{m,M}(g,g)\) as follows: 
\begin{equation}\label{Lp:QM}
\|Q_{m,M}(g,g)\|_{L^p_k} \le C_M \|g\|_{L^1_k}\|g\|_{L^p_k}, 
\quad 1\le p\le \infty
\end{equation}
(there will be no loss of moments since the 
kernel \(B_{m,M}\) is bounded). Therefore, 
\[
h \in{L^\infty([0,T], L^1_2 \cap L^2(\R^\di))}, 
\]
as soon as \(g\) is in the same space. 
The unique weak solution \(f\in L^\infty([0,T],L^1_2\cap L^2(\R^\di))\)  
of \eqref{fp} is then obtained from the following integral 
representation:  
\begin{equation}\label{fp:int}
f(v,t) = f_0(v) * E(v,t) + \int_0^{t} h(v,\tau) * E(v,t-\tau) \, d\tau, 
\end{equation}
where \(*\) denotes the convolution in \(v\), and 
\(E(v,t)\) is the fundamental solution of \eqref{fp}: 
\[
E(v,t)=\frac{1}{(4\pi t)^{\di/2}} e^{-\frac{|v|^2}{4t}-Mt}.  
\] 
The \(H^2\) regularity of \(f\) is then guaranteed by the classical 
parabolic regularity result \cite[Section 3.3]{LaUrSo}, 
and we have the bound
\begin{equation}
\label{fp:H2}
\|f\|_{H^2([0,T]\times\R^\di)} \le C_M 
(\|h\|_{L^2([0,T]\times\R^\di)} + \|f_0\|_{H^1(\R^\di)}). 
\end{equation}

We denote by \(\mathcal{T}\) the operator that maps \(g\) 
into \(f\). We next establish that for a certain choice 
of constants \(A_1\) and \(A_2\) this operator maps the set 
\begin{equation}
\label{def:B}
\begin{split}
B=\Big\{f & \in L^1([0,T]\times\R^\di)) \;\Big|\; 
f \ge 0, \;
\ir f \, dv=1, \; \ir f\,v\,dv = 0, 
\\
& \ir f\,|v|^2\,dv\le A_1,\;   \ir f^2\, dv \le A_2^2, 
\text{ for } a.\,a. \; t\in[0,T]\Big\}
\end{split}
\end{equation}
into itself. Indeed, the nonnegativity of \(f\) is evident from 
the integral representation \eqref{fp:int}, since \(h\ge0\). 
The mass and momentum normalization conditions follow easily, 
since for \(g\in B\) the collision term \(Q_{m,M}(g,g)\) integrates 
to zero when multiplied by \(1\) or \(v\). It remains to verify 
the last two conditions in \eqref{def:B}. 

For the first of these conditions, multiplying the equation 
\eqref{fp} by \(|v|^2\) and integrating by parts we obtain:
\begin{equation}
\label{fp:L12}
\begin{split}
\frac{d}{dt}&{\ir}f\,|v|^2\,dv \,+\, M \ir{f}\,|v|^2\,dv 
\le 2\di + M \ir{g}\,|v|^2\,dv 
\\
-\,& k \ir\ir {g}\,{g_*}\,|v|^2|v-v_*|_{m,M}\,dv\,dv_*
\le 2\di + (M-mk) \ir{g}\,|v|^2\,dv,
\end{split}
\end{equation}
where \(k=\epsilon_\di(1-\alpha^2)/4\). 
Therefore, taking \(g\) so that  
\[
\ir{g}\,|v|^2\,dv \le A'_1=\frac{2\di}{mk},
\]
yields the differential inequality
\[
\frac{d}{dt}{\ir}f\,|v|^2\,dv + M \ir{f}\,|v|^2\,dv 
\le M A'_1. 
\]
Then, by Gronwall's lemma,
\[
{\ir}f\,|v|^2\,dv \le 
\max\Big(A'_1,{\ir}f_0\,|v|^2\,dv\Big)
\]
Therefore, setting 
\[
A_1=\max\Big(A'_1,{\ir}f_0\,|v|^2\,dv\Big),
\]
we obtain the required estimate.  

To obtain a bound of \(f\) in \(L^2\) we integrate 
the equation against \(f\) and use the inequality 
\eqref{Lp:QM} to estimate \(Q_{m,M}(g,g)\): 
\begin{equation}\label{L2:fp}
\begin{split}
\frac{1}{2}\frac{d}{dt}\ir f^2 \, dv 
+ \ir |\nabla{f}|^2 \,dv
& + M \ir f^2 \, dv
\\
& \le {C_M}\|g\|_{L^1}\|g\|_{L^2}\|f\|_{L^2} + M\|g\|_{L^2}\|f\|_{L^2}.
\end{split}
\end{equation}
By Sobolev's embedding and interpolation, 
\[
\|\nabla{f}\|_{L^2} \ge K \|f\|_{L^{2*}} 
\ge K\|f\|^{1/\lambda}_{L^2}\|f\|^{-(1-\lambda)/\lambda}_{L^1}
\]
where \(0<\lambda<1\) is as in \eqref{interp:lambda}. 
Therefore, dividing \eqref{L2:fp} by \(\|f\|_{L^2}\) 
and taking into account that 
\(\|g\|_{L^1}=\|f\|_{L^1}=1\), we get  
\[
\frac{d}{dt}\|f\|_{L^2}
+K^2\|f\|^{2/\lambda-1}_{L^2}+M\|f\|_{L^2}
\le (C_M+M)\|g\|_{L^2}.
\]
Using the inequality
\[
K^2x^{2/\lambda-1} \ge \frac{1}{\eps}\,x-K_\eps,
\]
true for all \(\eps>0\), we find 
\begin{equation*}
\begin{split}
\frac{d}{dt}\|f\|_{L^2}+({1}/{\eps}+M)\|f\|_{L^2}
\le K_\eps + (C_M+M)\|g\|_{L^2}.
\end{split}
\end{equation*}
We then get by Gronwall's lemma: 
\begin{equation}\label{fpineq2}
\begin{split}
\|f\|_{L^2} \le \max \big( \|f_0\|_{L^2}, \beta 
+ \gamma\|g\|_{L^2} \big). 
\end{split}
\end{equation}
where
\begin{equation}\label{fpineq3}
\begin{split}
\beta=\frac{K_\eps}{1/\eps+M}
\quad\mathrm{and}\quad
\gamma=\frac{C_M+M}{1/\eps+M}.
\end{split}
\end{equation}
Choosing \(\eps<1/C_M\) we get \(\gamma<1\). Therefore, we obtain 
the inequality \(\|f\|_{L^2} \le A_2\) if we set
\[
A_2 = \max \big( \|f_0\|_{L^2}, {\beta}/(1-\gamma) \big). 
\]

It is straightforward to verify 
that the set \(B\) is convex and closed in the strong topology of 
\(L^1([0,T]\times\R^\di)\), using Fatou's lemma and the fact 
that the second moment in \(|v|\) 
is uniformly bounded for \(g\in B\). 
Further, the uniform in time bounds assumed in the definition 
of \(B\) imply the continuity of \(Q_{m,M}(g,g)\) in \(L^1\). 
We can then deduce easily that the solution operator 
\(\mathcal{T}\) itself is continuous, based on the representation 
\eqref{fp:int}. 
Finally, the bound for the second moment and the 
regularity estimate \eqref{fp:H2} imply that the operator
\(\mathcal{T}\) maps \(B\) into its compact subset. 
By the Schauder theorem, this proves the  existence of a fixed 
point for \(\mathcal{T}\) in \(B\), which is thereby a weak 
solution \(f_{m,M}\in L^\infty([0,T],L^1_2\cap L^2(\R^\di))\) 
of \eqref{fp:trunc}.

Our next goal is to pass to the limit as \(M\to\infty\) and then 
as \(m\to0\), to recover the solutions with the ``hard sphere''
collision kernel. To this end, we { will} show that the bounds
set forth in the apriori estimates hold for the fixed point 
solutions, and are uniform in \(M\) (and \(m\)).  First of all, 
using the computation \eqref{fp:L12} it is easy to conclude 
that the second moment is bounded uniformly in \(M\), as soon
as \(m>0\). Indeed, we obtain the following inequality 
for \(f=f_{m,M}\),
\begin{equation}
\label{tdep-ke:m}
\frac{d}{dt}\ir f |v|^{2}\,dv \le 2N - km \ir f |v|^2\, dv,
\end{equation}
so the required bound follows by Gronwall's lemma. 

Further, we see that for every \(m>0\), \(M>0\) and for every \(T>0\), 
the solutions are in \(L^\infty([0,T],L^1_{2p}(\R^\di))\), 
for every \(p>1\).  To see this we take \(K>0\) and introduce 
the truncated function 
\[
\Psi_{p,K}(x) = 
\left\{
\begin{array}{l}
x^p, \quad 0\le x < K
\\
K^p + p\, K^{p-1} (x-K), \quad x \ge K. 
\end{array}
\right.
\]
Then \(\Psi_{p,K}(x)\) is convex in \(x\), continuously 
differentiable, and has a bounded second derivative. 
It also verifies conditions 
\eqref{cond:psi-1}--\eqref{cond:psi-4}, so Lemma \ref{povz:conv}
applies. 
Taking \(\Psi_{p,K}(|v|^2)\) as a test function in the 
weak form of \eqref{fp:trunc} and arguing as in Lemma 
\ref{povz:pow}, we get
\begin{equation}
\label{weak-povz:K}
\begin{split}
\frac{d}{dt}\ir f\,\Psi_{p,K}(|v|^2)
+p(p-1)k_p m
\int\int_{\{|v|^2+|v_*|^2\le K\}} f f_*\, |v|^{2p}\,dv\,dv_*
\\
\le 2 pAM \ir f\,|v|^2\,dv \ir f \,|v|^{2p-2}\,dv
+2p(2p-2+\di)\ir f|v|^{2p-2}\,dv. 
\end{split}
\end{equation}
Therefore, if we take \(1<p\le2\), we can pass to the limit as \(K\to\infty\)
in \eqref{weak-povz:K} using the monotonicity with respect to \(K\) and  
the bound of \(f\) in \(L^1_2\). This implies
\[
f\in L^\infty([0,T],L^1_{2p}(\R^\di)), 
\]
for every \(1<p\le2\), with bounds generally dependent on \(m\) 
and \(M\).  By induction, the same property 
is extended to every \(p\ge0\). 

We see further that the bounds in \(L^1_{2p}\) are in fact 
independent on \(M\). Indeed, estimating the middle term 
in \eqref{weak-povz:K} using the inequality
\[
|v-v_*|_{m,M} \le  m+|v-v_*|
\] 
and following the arguments of Lemma \ref{povz:pow} we get
\begin{equation}
\label{weak-povz:m}
\begin{split}
\frac{d}{dt}\ir \!f\,|v|^{2p}\,dv
+p(p-1)k_p m \!\ir \!f |v|^{2p}\,dv
\le A_p\!\ir \!f\,|v|\,dv \ir \!f \,|v|^{2p}\,dv 
\\
+ \Big(pm\ir f\,|v|^2\,dv + 2p(2p-2+\di)\Big) 
\ir f \,|v|^{2p-2}\,dv . 
\end{split}
\end{equation}
This implies that for every \(T>0\) fixed and every \(p\ge0\), 
the bounds of \(f=f_{m,M}\) in 
\(L^\infty([0,T],L^1_{2p}(\R^\di))\) are 
independent of \(M\). 

Using the established \(L^1_{2p}\) bounds and the fact that 
\(f\in H^2([0,T]\times\R^\di)\) we can make rigorous the arguments
of Lemma \ref{reg:tdep-1} and then proceed as in Lemmas 
\ref{reg:tdep-2}-\ref{reg:tdep-3} obtaining that 
\[
f\in L^\infty([0,T],H^n_{2p}(\R^\di)),
\]
for every \(n=1,2...\), and every \(p\ge0\), with bounds 
independent on \(M\).  This will allow us to pass to the 
limit as \(M\to\infty\) in the weak form and to show that
the limit solutions satisfy the equation with the kernel 
\[
(m+|u|)\,b(u,\sigma).
\]
We can then substitute the computation \eqref{weak-povz:m}
by the argument of Lemma \ref{povz:pow} and find 
the bounds in \(L^\infty([0,T],L^1_{2p}(\R^\di))\) that are
independent on \(m\) and \(T\). Arguing as above we can 
then pass to the limit as \(m\to0\). The limit solution 
obtained in this step will then satisfy the equation with
the ``hard sphere'' kernel.

Finally, in order to treat the problem with the initial data 
\(f_0\in L^1_{2}\cap L\log{L}(\R^\di)\) we can take a sequence 
\(f^n\in C^\infty_0(\R^\di)\) that converges to \(f_0\) in 
\(L^1(\R^\di)\). Then, since the constants in the bounds 
for the energy and entropy from Section \ref{apr:basic}
are independent of \(n\), we can pass to the weak \(L^1\)-limit 
in the equations. 
The fact that the bounds of the solutions 
are independent of \(T\) allows us to continue the obtained 
solutions to \([0,2T]\), and by induction, to \([0,\infty)\).

To study the regularity of solutions with \(L^2\) initial 
data we use the parabolic regularity of the equation \cite{LaUrSo}
to find that \(f\in H^1([0,T]\times\R^\di)\), for any \(T>0\). 
Using this fact in combination with the bound 
\(f\in L^\infty([0,\infty),L^1_\rbeta(\R^\di))\) we can make 
rigorous the argument of Lemma \ref{reg:tdep-1} and then 
proceed as in Lemmas \ref{reg:tdep-2} and \ref{reg:tdep-3} 
to find the infinite differentiability of the solutions. 
\end{proof} 

We now turn our attention to the steady equation \eqref{be:stat}
and give a proof of Theorem \ref{exist:stat}. One of the possible 
approaches consists of adapting the arguments developed above for 
the time-dependent case. In fact, as a careful reader will easily 
check, practically all arguments in the above proof apply to the steady 
equation: the Gronwall lemma arguments will be replaced by the 
inequalities obtained by dropping the time-derivative terms.  
The only point that would need more careful attention is the 
moment estimate \eqref{weak-povz:m}, which is {\it not} uniform 
in \(T\). It can be replaced by a more elaborate argument 
for the moment bounds in the case of the truncated collision 
kernel. We will, however take another approach, which will 
allow us to obtain the existence of the steady problem as a 
consequence of the regularization properties of the time-dependent 
equation. 

\begin{proof}[Proof of Theorem \ref{exist:stat}]
The proof of Theorem \ref{exist:tdep} enables us to construct 
a semigroup on the convex set $C$ made of those functions in 
$L^1_2\cap L^2({\Bbb R}^N)$ with unit mass and zero mean. 
Denote it by $(S_t)_{t\geq 0}$. Our bounds imply that for all
$t>0$, the range of $S_t$ is compact in $C$. Therefore, for
all $n$ the equation
$$ f_n = S_{2^{-n}} f_n$$
is solvable by Schauder's theorem.
Since $f_n = S_1 f_n$, the sequence $f_n$ is contained in a fixed compact
of $C$, namely $S_1(C)$. We can therefore extract a subsequence
which converges towards some $f$. Now for all $k\leq n$ we have
$$ f_n = S_{2^{-k}}f_n $$
(because $2^{-k}$ is a multiple of $2^{-n}$), and we can pass
to the limit as $n\to\infty$ using the continuity of the semigroup, 
thereby obtaining 
$$ f = S_{2^{-k}} f, \qquad \text{for all }\quad k\geq 0. $$
Therefore $f = S_t f$ for all $t$ which is a sum of inverse
powers of~2. Since the set of such times forms a dense subset
of ${\mathbb R}_+$ and since the semigroup is continuous with
respect to $t$, we conclude that
$$ f = S_t f, \qquad \text{for all }\quad t\geq 0.$$
This ends the proof.
\end{proof}

\section{Uniqueness by Gronwall's lemma}

We next show that under the assumption that the initial 
data has the moment of order \(3\) finite, the solution 
to the time-dependent problem is unique. The proof uses 
an argument based on a certain cancellation property 
of the collision operator multiplied by \(\sgn(f)\) 
\cite{Ar2,Di} (see also \cite{MiWe} for discussion).  
We show that this property yields the 
desired result for the operator with inelastic collisions 
as well.  

\begin{theorem}
\label{thm:uniq}
Assume that \(f_0\in L^1_3(\R^\di)\); then the equation 
\eqref{be:tdep} with the initial condition \(f(\cdot,0)=f_0\) 
has at most one solution. 
\end{theorem}
\begin{proof}
Assume that \(f\) and \(g\) are solutions of \eqref{be:tdep}, 
with the same initial data \(f_0\). Set \(h=f-g\) and 
\(H=f+g\). Then \(h\) satisfies the equation
\begin{equation}\label{be:uniq}
\partial_t{h} - \Delta_v{h} = \frac{1}{2}\,(Q(h,H)+Q(H,h)),
\end{equation}
with the homogeneous initial data. 
Now take a function \(\psi_\eps(x)\), a continuous 
approximation of \(sgn(x)\).  We can take 
\[
\psi_\eps(x) = 
\left\{
\begin{array}{ll}
-1,\quad& x\le -\eps
\\
{x}/{\eps},\quad &-\eps < x \le \eps
\\
1,\quad & x > \eps. 
\end{array}
\right.
\]
Multiplying equation \eqref{be:uniq}
by \(\psi_\eps(h)\,(1+|v|^2)\) and integrating by parts we get
\begin{equation*}
\begin{split}
\frac{d}{dt}\ir h\,\psi_\eps(h) (1+|v|^2) \, dv 
+ \frac{1}{2\eps}\int_{\{|h| \le \eps\}} |\nabla{h}|^2 (1+|v|^2) \, dv
\\
- 2 \di \ir \phi_\eps(h) \,dv 
= \frac{1}{2} \ir (Q(h,H)+Q(H,h)) \, \psi_\eps(h)\, (1+|v|^2) \,dv, 
\end{split}
\end{equation*}
where
\[
\phi_\eps(x)=\int_0^x \psi(t)\,dt=\left\{
\begin{array}{ll}
-x+{\eps}/{2},\quad& x\le -\eps
\\
{x^2}/{2\eps},\quad &-\eps < x \le \eps
\\
x-{\eps}/{2},\quad & x > \eps. 
\end{array}
\right.
\]
To estimate the right-hand side we can adapt the argument that is 
known to work in the case of the elastic Boltzmann equation 
(cf. \cite{Ar2,Di}). 
Passing to the weak form we get:
\begin{equation}\label{4int}
\begin{split}
& \ir (Q(h,H)+Q( H,h))\,\psi_\eps(h)\,(1+|v|^2) \,dv
\\
& = \ir\ir\is H_* \,\big\{  h \,\psi_\eps(h')\,(1+|v'|^2) 
+ h \,\psi_\eps(h'_*)\,(1+|v'_*|^2)
\\
& \qquad\quad - \, h \,\psi_\eps(h)\,(1+|v|^2) 
- h \,\psi_\eps(h_*)\,(1+|v_*|^2)\,\big\}\,
|u|\,b(u,\sigma)\,d\sigma\,dv\,dv_*.
\end{split}
\end{equation}
Since \(|v'|^2+|v'_*|^2 \le |v|^2+|v_*|^2\), we can 
estimate the integrals 
of the first two terms in the braces as follows: 
\[
(2+|v|^2+|v_*|^2-c_\di\frac{1-\alpha^2}{4}|v-v_*|^2) \,
\ir\ir |h| \,H_*\, (2+|v|^2+|v_*|^2) \,
|v-v_*|\,dv\,dv_*.
\]
Subtracting the third term in the integral \eqref{4int} and noticing that 
\[
\big|\,h\psi_\eps(h)-|h|\,\big| 
\le |h|\,\chi_\eps(h),
\]
where \(\chi_\eps(x)\) is the characteristic function 
of the interval \([-\eps,\eps]\), we obtain the estimate for 
the first three terms:
\[
\begin{split}
\ir\ir |h| \,H_*\, (1+|v_*|^2) \,|v-v_*|\,dv\,dv_*
\\
+\ir H_* \int_{\{|h|\le\eps\}} h \,(1+|v|^2) |v-v_*|\,dv\,dv_*.
\end{split}
\]
The fourth term in \eqref{4int} contributes with another integral 
like the first one above, so we finally get
\begin{equation*}
\begin{split}
\frac{d}{dt}\ir h\,\psi_\eps(h) (1+|v|^2) \, dv 
\\
\le 2 \di \ir \phi_\eps(h) \,dv 
+\ir\ir |h| \,H_*\, (1+|v_*|^2) \,|v-v_*|\,dv\,dv_*
\\
+\frac{1}{2}\ir H_* \int_{\{|h|\le\eps\}} h \,(1+|v|^2) |v-v_*|\,dv\,dv_*.
\end{split}
\end{equation*}
Passing to the limit as \(\eps\to0\), 
we find:
\[
\begin{split}
&\frac{d}{dt}\ir |h|\,(1+|v|^2)\, dv 
\\
&\le 2 \di \ir |h| \,dv + \ir |h| \,(1+|v|^2)^{1/2} \,dv 
\ir H \,(1+|v|^2)^{3/2} \,dv  
\\
&\le C \ir |h| \,(1+|v|^2) \,dv, 
\end{split}
\]
since \(H\) is assumed to be bounded in \(L^1_3(\R^\di)\).
Now since \(h(0,v)=0\) it follows by Gronwall's lemma that 
\(h(t,v)=0\) for all times. 
\end{proof}

\begin{remark}
The uniqueness result of Theorem \ref{thm:uniq} is 
most certainly suboptimal. We believe that the 
uniqueness could be obtained in the class of initial 
conditions with finite mass and energy, with 
no additional assumptions, similarly to the 
the classical Boltzmann equation \cite{MiWe}. 
The main technical obstacle for such a result 
is extending the Povzner inequalities in the case 
inelastic collisions to the class of slowly growing 
piecewise linear functions \(\psi\) studied in 
\cite{MiWe}. We believe that this can be overcome 
with a more careful 
analysis of the inelastic collision mechanism. 
\end{remark}

\section{Lower bounds with overpopulated high energy tails}

In this section we obtain pointwise lower estimates
of solutions to \eqref{be:tdep} and \eqref{be:stat}
showing that the behavior of the high-energy tails 
of solutions is controlled below by ``stretched Maxwellians''
\(A\exp(-a|v|^{3/2})\). The bounds are established
by using the comparison principle based on the parabolic
(elliptic) structure of the equations. 
The following proposition establishes the particular role 
played by the ``stretched Maxwellians'': they can be used
as barrier functions in the comparison principle.

\begin{proposition}\label{h.ineq}
Let  \(g(v)\) be a nonnegative function with finite 
mass \(\rho_0 = \ir g\,dv\) and moment of order one, 
\(\rho_1 = \ir g \,|v|\,dv\). Then for every \(r > 0 \),
and  every  \(K>0\), there is a constant \(a>0\) such that 
the function  
\begin{equation*}
\begin{split}
h(v) = K e^{- a |v|^{3/2}}, 
\end{split}
\end{equation*}
satisfies
\begin{equation}
\label{ine1}
\begin{split}
\Delta {h} - Q^-(g,h) \ge 0, \quad \text{for all} \quad |v| > r,
\end{split}
\end{equation}
Further, choosing \(b>0\) large enough, the function 
\begin{equation*}
\begin{split}
h(v,t) = K e^{ - bt - a (1 +|v|^2)^{3/4}}. 
\end{split}
\end{equation*}
satisfies
\begin{equation}
\label{tine1}
\begin{split}
- \partial_t h + \Delta {h}   - Q^-(g,h)  \ge 0,
\end{split}
\end{equation}
for all \(t>0\) and all \(v\in\R^\di\). 
\end{proposition}
\begin{proof}
To prove inequality \eqref{ine1} we fix an \(r>0\), 
compute, 
\begin{equation*}
\begin{split}
\Delta {h} 
= \Big(\, \frac{9}{4}\, a^2 |v| 
- \frac{3(2\di-1)}{4}\, a |v|^{-1/2} \Big) \, h  
\end{split}
\end{equation*}
and use the estimate
\[
Q^-(g,h) = h\,(g * |v|) \le (\rho_1+\rho_0|v|)\, h,
\]
to obtain
\begin{equation}\label{ine2}
\begin{split}
\Delta {h} - Q(g,h)
&\ge \Big( \Big(\frac{9}{4}\, a^2-\rho_0) |v| - \rho_1 
- \frac{3(2\di-1)}{4} a |v|^{-1/2} \Big) \, h\,.
\end{split}
\end{equation}
If \(\frac{9}{4} \,a^2 \ge \rho_0\), 
the factor on right-hand side 
of \eqref{ine2} attains its minimum for \(|v|=r\).   
Therefore, inequality \eqref{ine1} holds for every 
\(a\ge a^*\), where \(a^*\) is the positive root of the 
quadratic equation
\[
\frac{9r}{4} \,a^2 
- \frac{3(2\di-1) r^{-1/2} }{4} \,a 
- (\rho_0 r +\rho_1) = 0
\]
For the time-dependent operator, denoting by 
$\jap{v}=(1 + |v|^2)^{\frac12}$,  we obtain
\begin{equation}\label{tine2}
\begin{split}
\Delta {h} - \partial_t h 
= \Big( \,\frac{3 a}{4}\, |v|^2 \jap{v}^{-1}
\big(3a + \jap{v}^{-1/2} \big) 
- \frac{3 \di a}{2}\jap{v}^{-1/2} + b  \Big) \, h.  
\end{split}
\end{equation}
Choosing \(a\) so that \(\frac{9}{4}\,a^2\ge \rho_0\) and then 
\(b \ge \frac{3\di a }{2}+\rho_0+\rho_1\)
we obtain inequality \eqref{tine1} and complete the proof.   
\end{proof}
The established property of the function \(h(v)\) is used 
in next lemma to obtain a comparison result for the 
steady equation.  
\begin{lemma}\label{mp:stat}
Let \(f \in L^1_1(\R^\di)\) be a nonnegative smooth solution 
to \eqref{be:stat}, with the mass \(\rho_0>0\).   
Then, there is a constant \(K>0\), such that 
\begin{equation}\label{compare}
f(v) \ge K e^{-2a |v|^{3/2}} , 
\end{equation}
for all \(v\in\R^\di\),
where \(a\) is a constant as in Proposition \ref{h.ineq}. 
\end{lemma}
\begin{proof}
Assuming the smoothness of the solution to \eqref{be:stat},
there is a constant \(c_0>0\) and a ball \(B(v_0,r_0)\) with 
\(v_0 \in \R^\di\) and \(r_0>0\), such that 
\begin{equation}
\label{in.ball}
f(v) \ge c_0>0, \quad \text{if}\quad v \in B(v_0,r_0). 
\end{equation}
The value of \(c_0\) (as well as \(r_0\) and \(v_0\)) 
depend on the solution \(f\) and use the fact that 
\(\rho_0>0\). 

Since equation \eqref{be:stat} is translation invariant, 
we can take \(g(v) = f(v+v_0)\); then
\begin{equation}\label{compare1}
\Delta g - Q^-(g, g) = \Delta g - (g* |v|) g \le 0.
\end{equation}
Applying Proposition \ref{h.ineq} to the function \(g(v)\)
with \(r=r_0\) we find the barrier function 
\(h(v)=c_0 \, \exp(-a|v|^{3/2})\), for which we have
\begin{equation}
\label{compare2}
\Delta h - Q^-(g,h) = \Delta h - ( g* |v| ) h  \ge  0, 
\quad \text{for } \;|v| >r_0.
\end{equation}
and
\begin{equation*}
g(v) \ge h(v), \quad \text{for} 
\;\; |v| \le r_0. 
\end{equation*}

Therefore, letting \(U(v)=g(v)-h(v)\), subtracting \eqref{compare2} from 
\eqref{compare1} we obtain 
the inequality
\begin{equation*}
\Delta U - (g* |v|) \,U \le 0, \quad |v|>r_0,
\end{equation*}
To prove that \(U(v) \ge 0\) everywhere we  
apply a form of a strong maximum principle (see, for example,  
\cite{GiTr}) to the operator
\[
\mathcal{L}\, U = \Delta U - \nu(U+h) \,U. 
\] 
We can reduce the problem to proving that \(U \ge 0\) in a bounded 
domain. Indeed, the decay conditions on \(f\) imply that for every 
\(\eps>0\) we can find \(R>0\) such that \(|U(v)|<\eps\) 
if \(|v| \ge R\). Then we have 
\[
\mathcal{L}(U+\eps) = \mathcal{L} \,U - \eps \nu(g) \le 0 , 
\quad r_0 < |v| < R
\] 
and \(U+\eps>0\) for \(|v|=r_0\) and \(|v|=R\). The strong 
maximum principle then implies that \(U+\eps \ge 0\) for all 
\(r_0 \le |v| \le R\). Letting \(\eps\) go to zero we get 
\[
U \ge 0, \quad\text{for all}\quad |v| \ge r_0. 
\] 
In view of the inequality \eqref{in.ball} this implies
\[
g(v) \ge c_0 \, e^{- a |v|^{3/2}},
\]
or, applied to the function \(f(v)\), 
\[
f(v) \ge c_0\, e^{- a |v-v_0|^{3/2}} 
\ge K e^{-2a|v|^{3/2}},
\]
with \(K=c_0 \, e^{-a |v_0|^{3/2}}\). 
This completes the proof of the lemma. 
\end{proof}

By using a version of the maximum principle for the 
parabolic operator, we obtain, in a similar fashion,
the pointwise lower bound for the time dependent problem.

\begin{lemma}\label{lower:comp:time}
Let \(f\in L^\infty([0,\infty),L^1_2(\R^\di))\) be 
a nonnegative smooth solution to \eqref{be:tdep}
with the initial data \(f_0 \ge c_0\exp(-a_0|v|^{3/2})\).   
Then, there are  positive constants \(K\), \(b\) and 
\(a\), generally depending on the solution, such that  
\[
f(v,t) \ge K e^{- bt -a |v|^{3/2}},
\]
for all \(t>0\) and all \(v\in\R^\di\). 
Further, if there is a constant \(c_1\) 
and  a ball \(B(v_0,r_0)\), such that   
\begin{equation*}
f(v,t) \ge c_1, \quad \text{if}\quad v \in B(v_0,r_0),
\end{equation*}
for all \(t\), then the lower bound 
\[
f(v,t) \ge K e^{-a |v|^{3/2}},
\]
holds uniformly in time, where now  \(K>0\), \(a>0\) and 
\(b>0\) will depend on \(c_1\), \(v_0\) and \(r_0\).
\end{lemma}
\begin{proof} To prove the first statement of the lemma we 
use the second part of Proposition \ref{h.ineq} and 
repeat the comparison arguments of Lemma \ref{mp:stat} taking 
\[
h(v,t) = K e^{-bt -a\jap{v}^{3/2}}
\] 
and using the strong maximum principle 
for the parabolic operator  on $U=f(v+v_0) - h(v)$
\[
\mathcal{L}\, U = \Delta U - \nu(f) \,U - \partial_t U. 
\] 
For the second part, the additional assumption made on 
\(f\) allows us to repeat the proof of Lemma \ref{mp:stat}
using the function \(h\) from \eqref{tine2}. 
\end{proof}

\begin{remark} 
It is tempting to conjecture that solutions to \eqref{be:stat} 
should satisfy a pointwise upper bound of the type 
\(K'\exp(-a'|v|^{3/2})\), for certain values of \(a'\) 
and \(K'\). However, the application of an argument 
based on the maximum principle requires estimating 
\(Q^+(f,f)\) pointwise, which is generally a difficult 
problem. Assuming a ``no-cancellation'' property 
in the spirit of the argument \cite{VaEr}, 
\begin{equation}
\label{big-assumpt}
Q^+(f,f) - Q^-(f,f) \le - k_\alpha Q^-(f,f),
\end{equation}
where \(k_\alpha>0\), a pointwise upper bound is indeed 
obtained by the maximum principle techniques. However, 
a justification of \eqref{big-assumpt} at the present time 
seems to be out of reach. Notice that quite recently
Bobylev {\it et al.} \cite{BoGaPa} were able to prove 
an upper bound ``in the \(L^1\) sense'', namely, that 
for a certain choice of \(a'>0\)
\[
\ir f(v)\exp(a'|v|^{3/2})\,dv < + \infty, 
\] 
which could possibly be a hint in favor of the pointwise 
bound hypothesis.
\end{remark}

\section{Concluding remarks}

We studied the existence, uniqueness and regularity for the
time-dependent equation \eqref{be:tdep} and the existence,
regularity for the steady equation \eqref{be:stat}. 
An important problem that remained beyond the scope of our 
study is the convergence of the time-dependent solutions to 
the steady ones as time approaches infinity. In fact, this 
remains a serious open problem, since {no Lyapunov functional
for the time-evolution is known to exist}.
A number of other interesting questions can be raised 
in connection to the obtained results. Are the steady
states unique up to a normalization? Do the steady 
solutions necessarily have radial symmetry? (This can 
be expected from the rotation invariance of the 
equations; the existence and regularity of radial 
solutions can be obtained by applying our analysis 
to the reduced one-dimensional problem, {as in} \cite{CeIlSt},
{ or just by working in spaces of radially symmetric functions}). 

We hope that the methods developed in the present 
work for the case of diffusion forcing could be useful 
for studying other problems involving the Boltzmann 
(Enskog) collision 
terms with other collision and driving mechanisms. In
particular, a generalization to the case of a heat 
bath including a friction term seems to be rather 
straightforward. (The lower bounds in that case are 
expected to be Maxwellians.) It is also likely that
applying the techniques of this paper should yield 
results for problems with the normal restitution 
coefficient dependent on the 
relative velocity \cite{BiShSwSw,BrPo}, which would 
allow us to study a broader range of physical phenomena. 

Another problem worth studying is the (quasi-)elastic 
limit \(\alpha\to1\). The steady states for the 
Boltzmann equation with elastic interactions (\(\alpha=1\)) 
and vanishing diffusion (\(\mu=0\)) are Maxwellians, 
while for every \(\mu>0\) and every \(\alpha<1\) we 
have a ``3/2'' lower bound. Obtaining quantitative 
information on the transition to the Maxwellian 
steady states would be valuable. We hope to address 
some of these questions in our future work.   

\section*{Acknowledgements}

A number of people have contributed in a very important way 
to the development of 
the present paper. We would like to thank C. Bizon, S. J. Moon, 
J. Swift and H. Swinney for discussions concerning the physical 
aspects of the problem. 
We are also thankful to A. V. Bobylev, J. A. Carrillo, 
C. Cercignani and S. Rjasanow for fruitful 
discussions and a number of suggestions that helped us 
to improve the presentation. The first author has been supported by  NSF 
under grant DMS 9971779 and by TARP 
under grant 003658-0459-1999; {the third author has been supported
by the HYKE European network, contract HPRN-CT-2002-00282.} 
Support from the  Texas Institute for Computational and Applied
Mathematics/Austin is also gratefuly acknowledged. 

\bibliographystyle{hsiam}
\bibliography{gran12a}

\end{document}